\documentclass[]{siamart220329}
\usepackage[T1]{fontenc}
\usepackage{xr}
\usepackage[utf8]{inputenc}
\usepackage{amsfonts}
\usepackage{amsmath,amssymb}
\usepackage[caption=false]{subfig} % <- Preamble
\usepackage{tikz,pgfplots}
\usepackage{multirow}
\usepackage{makecell}
\usepackage{tikz-3dplot}
\usetikzlibrary{math}
\usepackage{bm}
\usepackage{rotating}
\usepackage{algorithm}
\usepackage{algpseudocode}
\usepackage{graphicx,epstopdf} % <- Preamble
\usepackage{enumerate}
\usepackage{enumitem}
\usepackage{graphicx}

\usepackage{hyperref}
\usepackage{relsize}
\usepackage{adjustbox}
\usepackage{cleveref}
\usepackage{tcolorbox}
\usepackage{dsfont}
\usepackage{accents}

\usepackage{multirow}
\usepackage{booktabs}
\usepackage{graphicx}

\usepackage{pdflscape}

\usepackage{perpage}
\MakePerPage{footnote}

\usepackage{colortbl}
\usetikzlibrary{fadings}
\ifpdf
  \DeclareGraphicsExtensions{.eps,.pdf,.png,.jpg}
\else
  \DeclareGraphicsExtensions{.eps}
\fi

\headers{New Hierarchical Matrices for fast kernel MVP with Applications}{Ritesh Khan, Sivaram Ambikasaran}

\title{New hierarchical matrices for fast kernel matrix-vector products with Applications to Gaussian processes and kernel ridge regression}

\author{Ritesh Khan\thanks{Department of Numerical Mathematics, Faculty of Mathematics and Physics, Charles University, Prague, Czech Republic. \email{ritesh.khan@matfyz.cuni.cz}, \email{khanritesh28@gmail.com}}
\and Sivaram Ambikasaran\thanks{
Department of Data Science \& Artificial Intelligence, Wadhwani School of Data Science \& Artificial Intelligence,
Department of Mathematics and  IIT Madras, Chennai, India. \email{\lowercase{sivaambi@alumni.stanford.edu}}}}

\newsiamremark{remark}{Remark}

\hypersetup{
    colorlinks=true,
    urlcolor=blue,
    citecolor = magenta
    }
\usetikzlibrary{arrows} 
\usetikzlibrary{decorations.markings}

\pgfplotsset{compat=1.16} 

\newcommand{\dsum}{\displaystyle\sum}

\newcommand{\xb}{\pmb{x}}

\newcommand{\bkt}[1]{\left(#1\right)}
\newcommand{\abs}[1]{\left\lvert#1\right\rvert}
\newcommand{\magn}[1]{\left\lVert#1\right\rVert}

\newcommand{\ceil}[1]{\left\lceil#1\right\rceil}

\newcommand{\Rb}{\mathbb{R}}
\newcommand{\Cb}{\mathbb{C}}

\newcommand{\cmt}[1]{\iffalse {#1} \fi}

\makeatletter
\newcommand*{\addFileDependency}[1]{% argument=file name and extension
  \typeout{(#1)}% latexmk will find this if $recorder=0 (however, in that case, it will ignore #1 if it is a .aux or .pdf file etc and it exists! if it doesn't exist, it will appear in the list of dependents regardless)
  \@addtofilelist{#1}% if you want it to appear in \listfiles, not really necessary and latexmk doesn't use this
  \IfFileExists{#1}{}{\typeout{No file #1.}}% latexmk will find this message if #1 doesn't exist (yet)
}
\makeatother

\newcommand{\mvp}{matrix-vector product }

\DeclareUnicodeCharacter{2212}{-}

\begin{document}

\maketitle

\begin{abstract}
This article presents two new multilevel hierarchical matrices to perform fast kernel matrix-vector product in $d$ dimensions $(d>1)$, namely $\mathcal{H}_{hw}^{\prime}$, $\mathcal{H}_{hw}^{\prime \prime}$. We introduce a novel \emph{hybrid weak admissibility} condition to construct these hierarchical matrices. More precisely, at coarser levels of the tree (i.e., larger blocks), we use a weak admissibility condition introduced in \cite{khan2022numerical}, wherein \textbf{far-field} and \textbf{vertex-sharing} clusters are considered as admissible, while at the finer levels (i.e., smaller blocks), the \textbf{adjacent} clusters are treated as admissible. The admissible blocks at coarser levels of $\mathcal{H}_{hw}^{\prime}$ and $\mathcal{H}_{hw}^{\prime \prime}$ are represented using \emph{semi-nested} and \emph{fully nested} bases, respectively. In both hierarchical matrices, the finer levels at which adjacent clusters are treated as admissible are constructed in a non-nested fashion. Thus, they can be viewed as a hybrid of $\mathcal{H}$ and $\mathcal{H}^2$ matrices. We use adaptive cross approximation (ACA) and nested cross approximation (NCA) to construct non-nested and nested bases, respectively. The fast matrix-vector products are performed using the proposed hierarchical matrix formats. To evaluate the performance and usefulness of the proposed hierarchical matrices, we conduct some numerical experiments. The numerical results show that they outperform their \emph{uniform admissibility} counterparts (weak admissibility condition of \cite{khan2022numerical} is applied at all levels) in terms of storage and matrix-vector product cost, demonstrating the effectiveness of the proposed hybrid formats. We also show that the proposed formats are competitive with the standard (or strong) admissibility-based $\mathcal{H}^2$-matrices. Furthermore, we showcase some interesting applications of the proposed hierarchical matrices to (i) \emph{fast} principal component analysis in Gaussian processes and (ii) \emph{fast} kernel ridge regression, enabled by \emph{fast} iterative solver and \emph{fast} Nystr\"om preconditioner.
\end{abstract}



\section{Introduction}
Kernel matrices arise in a variety of numerical and statistical applications, such as PDEs \cite{greengard1987fast,ho2013hierarchical, massei2022hierarchical, ambikasaran2025hybrid}, inverse problems \cite{harlim2020kernel, ambartsumyan2020hierarchical}, machine learning \cite{gray2000,ker_den,gp_book}, etc. Given $N$ points $\{\pmb{x}_i \}_{i=1}^N$, an $N \times N$ kernel matrix is defined by
$$K\bkt{i,j}=f \bkt{\pmb{x}_i,\pmb{x}_j}, \qquad 1\leq i,j\leq N, $$
where $f$ is a prescribed kernel function. Such matrices are typically dense, and consequently naive matrix computations, such as matrix-vector product and the solution of linear systems, become prohibitive for large $N$. However, kernel matrices possess rank-structuredness that can be exploited to accelerate matrix computations. Over the past few decades, substantial research \cite{greengard1987fast,hackbusch1999sparse,borm2003hierarchical,xia2010fast,gillman2012direct,amestoy2017complexity} has been devoted to developing fast algorithms that exploit this structure, among which hierarchical matrices are a popular framework. Hierarchical matrices ($\mathcal{H}$-matrices) are data-sparse approximations of dense matrices by partitioning the matrix into blocks and approximating selected off-diagonal blocks by low-rank representations, thereby reducing both storage and arithmetic complexity. Some commonly used hierarchical matrix formats in BEM and kernel methods include $\mathcal{H}$-matrices \cite{hackbusch1999sparse}, $\mathcal{H}^2$-matrices \cite{borm2009construction}, HODLR (Hierarchically Off-Diagonal Low-Rank) matrices \cite{ambikasaran2013mathcal,ambikasaran2019hodlrlib}, HSS (Hierarchically Semi-Separable) matrices \cite{xia2010fast,chandrasekaran2007fast}, and HBS (Hierarchically Block-Separable) matrices \cite{gillman2012direct}. Hierarchical matrix formats can be distinguished by the admissibility conditions (e.g., standard admissibility \cite{borm2003hierarchical} and weak admissibility \cite{weak_hackbusch2004}) used to identify low-rank blocks and by the construction of the corresponding bases (e.g., non-nested \cite{hackbusch2000sparse} and nested \cite{borm2009construction}). It is worth noting that the nested bases construction in $\mathcal{H}^2$-matrices offers better computational efficiency than the non-nested bases construction used in $\mathcal{H}$-matrices.

Most existing hierarchical matrix works for multidimensional/higher-dimensional problems employ the standard (or strong) admissibility condition, in which admissible clusters are well-separated or far-field. In our recent works \cite{khan2022numerical, nHODLRdD}, we developed $\mathcal{H}$ and $\mathcal{H}^2$-matrices based on a weak admissibility condition \cite{khan2022numerical} that additionally treats \textbf{vertex-sharing} clusters as admissible. Numerical experiments showed that these formats can achieve better computational efficiency than standard admissibility-based $\mathcal{H}$-matrices while maintaining comparable accuracy. However, one computational downside of these $\mathcal{H}$-matrix formats for multidimensional problems is that the number of admissible blocks (or translation operators) can grow rapidly with the tree depth, while many blocks (diagonal and neighbor blocks) remain dense at the leaf level.

In this article, we introduce a hybrid weak admissibility condition based on a prior switching level $\ell<L$, where $L$ denotes the tree depth. The weak admissibility condition \cite{khan2022numerical}, under which both far-field and vertex-sharing clusters are admissible, is applied up to level $\ell$ (coarser), while at the remaining (finer) levels, adjacent clusters are treated as admissible. This construction reduces the number of admissible blocks at the finer levels, and only diagonal blocks are dense at the leaf level $L$. Thus, lowering the associated computational and storage costs. We develop two different $\mathcal{H}$-matrix formats based on our hybrid weak admissibility condition, using nested and non-nested bases constructions. We also provide an explicit condition under which the proposed hybrid $\mathcal{H}$-matrix formats are advantageous over their uniform counterparts.

The following are the main contributions of this article.

\begin{enumerate}
    \item We introduce a novel hybrid weak admissibility condition, under which far-field and vertex-sharing cluster interactions are admissible at coarser levels of the hierarchy, while at finer levels adjacent clusters are treated as admissible.
    \item We propose a new hierarchical matrix format, $\mathcal{H}_{hw}^{\prime}$, for fast matrix-vector product based on the hybrid weak admissibility condition. At coarser levels of the hierarchy, nested bases are employed for far-field interactions, while vertex-sharing interactions are represented using non-nested bases. At finer levels, all admissible interactions are represented using non-nested bases.
    \item We further enhance the $\mathcal{H}_{hw}^{\prime}$ format by employing nested bases for the vertex-sharing interactions at coarser levels, resulting in the improved hierarchical matrix format $\mathcal{H}_{hw}^{\prime \prime}$.
    \item Numerical experiments show the advantages of the proposed hybrid hierarchical matrix formats over conventional hierarchical matrix formats based on the (uniform) weak admissibility and standard admissibility conditions. We also apply the proposed hierarchical matrix formats to accelerate principal component analysis in Gaussian process and kernel ridge regression problems. As part of this work, we also make the implementation publicly available at \url{https://github.com/riteshkhan/Hhweak}.
\end{enumerate}

The remainder of this article is organized as follows. \Cref{sec:prelim} reviews hierarchical matrices under different admissibility conditions and discusses low-rank approximation techniques. In \Cref{sec:hybrid_matrix}, we introduce a new hybrid weak admissibility condition and propose two hierarchical matrix formats based on it. \Cref{sec:num_results} presents numerical experiments demonstrating the advantages of the proposed hybrid formats over the uniform formats, as well as their application to Gaussian process and kernel ridge regression problems. Finally, \Cref{sec:conclusion} concludes the article.

\section{Hierarchical matrices} \label{sec:prelim}
This section discusses hierarchical matrices based on various admissibility conditions and bases for admissible blocks. It also presents an overview of low-rank approximation techniques, ACA and NCA.

\subsection{Tree construction} 
Let $\Omega \subset \mathbb{R}^d$ denote a domain containing $N$ particles (or points), whose interactions are represented by a dense matrix $K \in \mathbb{R}^{N \times N}$. We further assume that $\Omega$ is contained within a $d$-dimensional hypercube $C$. We subdivide the hypercube $C$ using a tree data structure to construct hierarchical representations. Depending upon the particle distribution or requirement, one can use different trees, such as balanced tree \cite{grasedyck2003construction,hmatrix_book}, unbalanced tree \cite{hmatrix_book}, $k$-d tree \cite{ambikasaran2019hodlrlib}, etc. \textbf{For simplicity, we consider the $2^d$ balanced tree (e.g., for $d=3$, it is balanced oct tree) throughout this article}. At level $0$ (the root level), the hypercube $C$ itself. Each hypercube $C_i^{\bkt{l}}$ at level $l$ is subdivided into $2^d$ finer hypercubes $C_{\alpha}^{\bkt{l+1}}$, where $\alpha = 1,2,\dots, 2^d$, which belong to level $l+1$ of the tree. The former is the parent of the latter, and the latter are its children. The cluster $\mathcal{C}_i^{\bkt{l}}$ is defined as the collection of the indices corresponding to the particles in the hypercube $C_i^{\bkt{l}}$. Clusters at the same level that share a common parent are called siblings. The set of all siblings of $C_i^{\bkt{l}}$ is denoted by $siblings (C_i^{\bkt{l}})$. We stop the subdivision at a level $L$ (the leaf level) of the tree when each finest hypercube contains at most $n_{max}$ particles, where $n_{max}$ is a user-specified threshold that defines the maximum number of particles at each finest hypercube. The finest hypercubes at level $L$ are called the leaves. Let $\mathcal{T}^L$ represent the balanced $2^d$ tree with depth $L$. \Cref{fig:tree_div} illustrates the balanced quad tree at different levels in $2$D. It can be noted that $N \leq n_{max} 2^{d L} \implies L =  \ceil{\log_{2^d} \bkt{N/n_{max}}}$.

\begin{tikzfadingfrompicture}[name=rndpts]
    \fill[transparent!0] foreach ~ in {1,...,800}{(rand,rand) circle (0.01)};
\end{tikzfadingfrompicture}
\begin{tikzfadingfrompicture}[name=pivots]
    \fill[transparent!0] foreach ~ in {1,...,250}{(rand,rand) circle (0.02)};
\end{tikzfadingfrompicture}
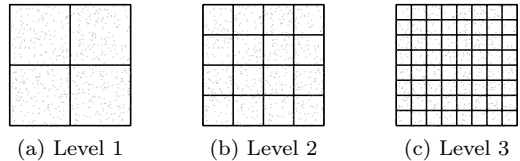
\begin{figure}[H]
    \centering
    % \subfloat[At level 0]{
    % \begin{tikzpicture}[scale=2]
    %     \node[text=red] at (0.5,0.5) {$0$};
    %     \draw[line width=0.2mm,  black] (0, 0) grid (1, 1);
    %     \draw (0,0) rectangle (1,1);
    %     \fill[black!20,scope fading=rndpts, fit fading=true] (0,0) rectangle (1,1);
    %     \node[text=red] at (0.5,0.5) {$0$};
    % \end{tikzpicture}
    % } \qquad
    \subfloat[Level 1]{
    \begin{tikzpicture}[scale=0.8]
        % \node[text=red] at (0.5,0.5) {$0$};
        % \node[text=red] at (1.5,0.5) {$1$};
        % \node[text=red] at (1.5,1.5) {$2$};
        % \node[text=red] at (0.5,1.5) {$3$};
        \draw[line width=0.2mm,  black] (0, 0) grid (2, 2);
        \draw (0,0) rectangle (2,2);
        \fill[black!20,scope fading=rndpts, fit fading=true] (0,0) rectangle (2,2);
        \draw[line width=0.2mm,  black] (0, 0) grid (2, 2);
        
        % \node[text=red] at (0.5,0.5) {$0$};
        % \node[text=red] at (1.5,0.5) {$1$};
        % \node[text=red] at (1.5,1.5) {$2$};
        % \node[text=red] at (0.5,1.5) {$3$};
    \end{tikzpicture}
    } \qquad
    \subfloat[Level 2]{
    \begin{tikzpicture}[scale=0.4]
        % \node[text=red] at (0.5,0.5) {$0$};
        % \node[text=red] at (1.5,0.5) {$1$};
        % \node[text=red] at (1.5,1.5) {$2$};
        % \node[text=red] at (0.5,1.5) {$3$};

        % \node[text=red] at (2.5,0.5) {$4$};
        % \node[text=red] at (3.5,0.5) {$5$};
        % \node[text=red] at (3.5,1.5) {$6$};
        % \node[text=red] at (2.5,1.5) {$7$};
        
        % \node[text=red] at (2.5,2.5) {$8$};
        % \node[text=red] at (3.5,2.5) {$9$};
        % \node[text=red] at (3.5,3.5) {$10$};
        % \node[text=red] at (2.5,3.5) {$11$};

        % \node[text=red] at (0.5,3.5) {$15$};
        % \node[text=red] at (1.5,3.5) {$14$};
        % \node[text=red] at (1.5,2.5) {$13$};
        % \node[text=red] at (0.5,2.5) {$12$};
        
        \draw[line width=0.2mm,  black] (0, 0) grid (4, 4);
        \draw (0,0) rectangle (4,4);
        \fill[black!20,scope fading=rndpts, fit fading=true] (0,0) rectangle (4,4);
        \draw[line width=0.2mm,  black] (0, 0) grid (4, 4);
        
        % \node[text=red] at (0.5,0.5) {$0$};
        % \node[text=red] at (1.5,0.5) {$1$};
        % \node[text=red] at (1.5,1.5) {$2$};
        % \node[text=red] at (0.5,1.5) {$3$};

        % \node[text=red] at (2.5,0.5) {$4$};
        % \node[text=red] at (3.5,0.5) {$5$};
        % \node[text=red] at (3.5,1.5) {$6$};
        % \node[text=red] at (2.5,1.5) {$7$};
        
        % \node[text=red] at (2.5,2.5) {$8$};
        % \node[text=red] at (3.5,2.5) {$9$};
        % \node[text=red] at (3.5,3.5) {$10$};
        % \node[text=red] at (2.5,3.5) {$11$};

        % \node[text=red] at (0.5,3.5) {$15$};
        % \node[text=red] at (1.5,3.5) {$14$};
        % \node[text=red] at (1.5,2.5) {$13$};
        % \node[text=red] at (0.5,2.5) {$12$};
    \end{tikzpicture}
    } \qquad
    \subfloat[Level 3]{
    \begin{tikzpicture}[scale=0.2]
        \draw[line width=0.2mm,  black] (0, 0) grid (8, 8);
        \draw (0,0) rectangle (8,8);
        \fill[black!20,scope fading=rndpts, fit fading=true] (0,0) rectangle (8,8);
        \draw[line width=0.2mm,  black] (0, 0) grid (8, 8);
    \end{tikzpicture}
    } \qquad 
    \caption{Hierarchical subdivision of the hypercube $C$ using balanced quad tree at different levels in $2$D $(d=2)$.}
    \label{fig:tree_div}
\end{figure}

\subsection{Admissibility conditions}
The low-rank block matrices of hierarchical matrices are identified using geometrical information known as admissibility conditions. Let $X$ and $Y$ be two clusters at the same level of a balanced $2^d$-tree, with corresponding hypercubes $C_i$ and $C_j$, respectively. For notational convenience, the superscript $``l"$ is occasionally omitted when both clusters are at the same level. If the admissibility condition is satisfied, the block matrix of $K$ corresponding to $X$ and $Y$, $K_{X,Y} \in \Cb^{\abs{I} \times \abs{J}}$, is approximated by a low-rank matrix with a prescribed tolerance $\epsilon$. Two commonly used admissibility conditions in the literature are the standard (or strong) admissibility condition and the weak admissibility condition. 

\subsubsection{Standard admissibility condition}
The block $X \times Y$, or the pair of clusters $\bkt{X,Y}$, is said to be admissible under the standard (or strong) admissibility condition, if the following condition holds:
\begin{align} \label{eq:H_admiss}
    \text{adm}_{\eta} \bkt{X, Y} = \text{true} \iff \dfrac{\min(\text{diam}(C_i), \text{diam}(C_j))}{\text{dist}(C_i,C_j)} \leq \eta, 
\end{align}
where 
\begin{align*}
    \text{diam}(U)= \sup\{\|\alpha- \beta\|_{2}:\alpha, \beta \in U \}
\end{align*}
and
\begin{align*} \label{eq:usual_distance}
    \text{dist}(U,V) = \inf\{\|u-v\|_{2}:u\in U, v\in V\}.
\end{align*}

The admissibility parameter $\eta>0$ controls the separation distance and the number of admissible pairs. From \cref{eq:H_admiss}, it follows that the smaller the $\eta$, the further apart the admissible clusters; conversely, the larger the $\eta$, the closer the admissible clusters.

If the admissibility condition is not satisfied, the corresponding cluster pair is inadmissible. The self-cluster and the neighboring clusters are inadmissible clusters.

The interaction list of a cluster $X$ is defined as
\begin{align*} \label{eq:IL_strong}
    \mathcal{IL}_{\eta} \bkt{X} := \big \{Y : \text{adm}_{\eta} \bkt{X, Y} = \text{true} \text{ and } \text{adm}_{\eta} \bkt{\text{parent} \bkt{X}, \text{parent} \bkt{Y}} = \text{false} \big \}.
\end{align*}

The neighbors list of $X$ is given by
\begin{align*} \label{eq:N_strong}
    \mathcal{N}_{\eta} \bkt{X} := \big \{ Y : \text{adm}_{\eta} \bkt{X, Y} = \text{false} \text{ and } Y \neq X \big \}.
\end{align*}

\begin{figure}[H]
\begin{center}
\subfloat[Level $1$]{
    \begin{tikzpicture}[scale=0.8]
        \fill [red] (1,0) rectangle (2,1);
        \fill [orange] (0,0) rectangle (1,1);
        \fill [orange] (1,1) rectangle (2,2);
        \fill [orange] (0,1) rectangle (1,2);
        \draw[step=1cm, black] (0, 0) grid (2, 2);
        % \node[anchor=north] at (1.5,0.7) {\tiny $\mathcal{C}^{(1)}$};
    \end{tikzpicture}
    \label{fmm_interaction_1}
    }\quad 
    \subfloat[Level $2$]{
    \begin{tikzpicture}[scale=0.4]
        \draw[black] (0, 0) grid (4, 4);
        \fill [green] (0,0) rectangle (1,1);
        \fill [green] (0,1) rectangle (1,2);
        \fill [green] (0,2) rectangle (1,3);
        \fill [green] (0,3) rectangle (1,4);
        \fill [green] (1,3) rectangle (2,4);
        \fill [green] (2,3) rectangle (3,4);
        \fill [green] (3,3) rectangle (4,4);
        \fill [red] (2,1) rectangle (3,2);
        \fill [orange] (1,0) rectangle (2,1);
        \fill [orange] (1,1) rectangle (2,2);
        \fill [orange] (1,2) rectangle (2,3);
        \fill [orange] (2,2) rectangle (3,3);
        \fill [orange] (3,2) rectangle (4,3);
        \fill [orange] (2,0) rectangle (3,1);
        \fill [orange] (3,0) rectangle (4,1);
        \fill [orange] (3,1) rectangle (4,2);
        % \node[anchor=north] at (2.5,1.8) {\tiny $\mathcal{C}^{(2)}$};
        \node[anchor=north] at (1.5,2.8) {};
        \node[anchor=north] at (1.5,1.8) {};
        \node[anchor=north] at (.5,1.8) {};
        \draw[black] (0, 0) grid (4, 4);
        \draw[line width=0.5mm,  black] (0, 0) rectangle (2, 2);
        \draw[line width=0.5mm,  black] (2, 2) rectangle (4, 4);
        \draw[line width=0.5mm,  black] (2, 0) rectangle (4, 2);
        \draw[line width=0.5mm,  black] (0, 2) rectangle (2, 4);
    \end{tikzpicture}
    \label{fmm_interaction_2}
    }\quad 
    \subfloat[Level $3$]{
    \begin{tikzpicture}[scale=0.4]
        % \draw[black] (0, 0) grid (4, 4);
        % \fill [green!25] (0,0) rectangle (1,1);
        % \fill [green!25] (0,1) rectangle (1,2);
        % \fill [green!25] (0,2) rectangle (1,3);
        % \fill [green!25] (0,3) rectangle (1,4);
        % \fill [green!25] (1,3) rectangle (2,4);
        % \fill [green!25] (2,3) rectangle (3,4);
        % \fill [green!25] (3,3) rectangle (4,4);
        % \node[anchor=north] at (2.28,2) {\tiny $\mathcal{C}^{(3)}$};
        \node[anchor=north] at (1.5,2.8) {};
        \node[anchor=north] at (1.5,1.8) {};
        \node[anchor=north] at (.5,1.8) {};
        \draw[black] (0, 0) grid (4, 4);
        \draw[line width=0.5mm,  fill=green] (1, 0) rectangle (2, 1);
        \draw[line width=0.5mm,  fill=green] (1, 1) rectangle (2, 2);
        \draw[line width=0.5mm,  fill=green] (2, 2) rectangle (3, 3);
        \draw[line width=0.5mm,  fill=green] (3, 2) rectangle (4, 3);
        \draw[line width=0.5mm,  fill=green] (2, 0) rectangle (3, 1);
        \draw[line width=0.5mm,  fill=green] (3, 0) rectangle (4, 1);
        \draw[line width=0.5mm,  fill=green] (3, 1) rectangle (4, 2);
        \draw[line width=0.5mm,  fill=green] (1, 2) rectangle (2, 3);
        \fill [orange] (1.5,1) rectangle (3,2.5);
        \fill [red] (2,1.5) rectangle (2.5,2);
        % \node[anchor=north] at (2.28,2.1) {\tiny $\mathcal{C}^{(3)}$};
        \draw[step=0.5cm, black] (1, 0) grid (4, 3);
        \draw[line width=0.5mm] (1, 0) rectangle (2, 1);
        \draw[line width=0.5mm] (1, 1) rectangle (2, 2);
        \draw[line width=0.5mm] (2, 2) rectangle (3, 3);
        \draw[line width=0.5mm] (3, 2) rectangle (4, 3);
        \draw[line width=0.5mm] (2, 0) rectangle (3, 1);
        \draw[line width=0.5mm] (3, 0) rectangle (4, 1);
        \draw[line width=0.5mm] (3, 1) rectangle (4, 2);
        \draw[line width=0.5mm] (1, 2) rectangle (2, 3);
    \end{tikzpicture}
    \label{fmm_interaction_3}
    }\quad
    \subfloat{
        \begin{tikzpicture}
            [%%%%%%%%%%%%%%%%%%%%%%%%%%%%%%
            box/.style={rectangle,draw=black, minimum size=0.2cm},scale=0.1
            ]%%%%%%%%%%%%%%%%%%%%%%%%%%%%%%
            \node[box,fill=red,,font=\tiny,label=right:Cluster considered,  anchor=west] at (8,0){};
            \node[box,fill=orange,,font=\tiny,label=right:Neighbor cluster,  anchor=west] at (8,4){};
            \node[box,fill=green,,font=\tiny,label=right:Far-field cluster,  anchor=west] at (8,8){};
            % \node[box,fill=green!50,,font=\tiny,label=right:Far-field cluster at parent level,  anchor=west] at (25,6){};
        \end{tikzpicture}
    }
    \caption{The interaction and neighbor lists of a cluster at different levels in $2$D.}
    \label{fig:h_interaction}
\end{center}
\end{figure}
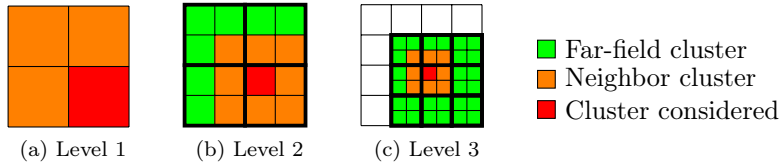

\Cref{fig:h_interaction} illustrates the neighbors and interaction lists of a cluster (red cluster) at different levels of the tree in $2$D, when we set $\eta = \sqrt{2}$ in \cref{eq:H_admiss}.

 The corresponding hierarchical matrix at different levels is shown in \Cref{fig:FMM2D_matrix}.

An upper bound on the cardinalities of the interaction list and the neighbor list of a cluster in a balanced $2^d$-tree can be derived \cite{weak_hackbusch2004}. The maximum sizes of the interaction list and the neighbor list for a cluster at \emph{any level} in $d$ dimensions are, respectively, as follows
\begin{align} \label{eq:H_matrix_IL_N}
   \abs{\mathcal{IL}_{\eta} \bkt{X}} \leq \bkt{2^d-1} \bkt{1 + {2\sqrt{d}}/{\eta}}^d, \qquad \abs{\mathcal{N}_{\eta} \bkt{X}} \leq \bkt{1 + {2\sqrt{d}}/{\eta}}^d - 1.
\end{align}

In particular, $\eta = \sqrt{d}$ in \cref{eq:H_admiss} yields a structure identical to balanced FMM (analytic counterpart of $\mathcal{H}^2$-matrices) \cite{greengard1987fast}, leading to

\begin{align} \label{eq:FMM_matrix_IL_N}
   \abs{\mathcal{IL}_{\sqrt{d}} \bkt{X}} \leq \bkt{2^d-1} 3^d, \qquad \abs{\mathcal{N}_{\sqrt{d}} \bkt{X}} \leq \bkt{3^d - 1}.
\end{align}

\subsubsection{Weak admissibility condition}
The notion of weak admissibility was introduced by Hackbusch et al. \cite{weak_hackbusch2004} for \emph{one-dimensional} problems, where they compress the sub-matrices corresponding to the adjacent clusters (intervals). The notion of weak admissibility doesn't look at the distance between the clusters; instead, it is defined as non-self or non-overlapping clusters being admissible. To be more precise, in terms of matrix setting, this represents all the off-diagonal sub-matrices as low-rank. Hierarchical matrix classes such as HODLR (non-nested bases) \cite{ambikasaran2013mathcal,ambikasaran2019hodlrlib} and HSS (nested bases)~\cite{xia2010fast,chandrasekaran2007fast} are based on Hackbusch's weak admissibility condition. 

Hackbusch et al. discuss the notion of weak admissibility only for one-dimensional problems. However, the extension of this notion to $d$ dimensions ($d>1$), or more generally to multi-dimensional settings, is not discussed in \cite{weak_hackbusch2004}. A straightforward extension to $d$ dimensions, i.e., compressing all the off-diagonal sub-matrices, will not result in a \emph{quasi-linear} complexity hierarchical matrix algorithm. This is because, in $d$ dimensions with $d>1$, except for the vertex-sharing interaction sub-matrix, the rank of all other adjacent interactions sub-matrices grows with some positive power of $N$ \cite{ho2013hierarchical,kandappan2022hodlr2d,khan2022numerical}. Thus, the notion of weak admissibility introduced by Hackbusch et al. is well-suited to one-dimensional problems but requires modification for multi-dimensional problems.

In our recent article \cite{khan2022numerical}, we show that the rank of not just the \emph{far-field} but also the \emph{vertex-sharing} interaction sub-matrices does not scale with any positive power of $N$. Therefore, the admissibility of \emph{far-field} and  \emph{vertex-sharing} clusters could be considered a way to extend the notion of weak admissibility in higher dimensions. We extend the notion of weak admissibility condition to any dimension $d$ as follows: 

\subsubsection*{Weak admissibility condition in $d$ dimensions}  \label{subsec:weadk_d_dim}
The block $X \times Y$, or the pair of clusters $\bkt{X,Y}$, is said to be admissible under our weak admissibility condition if the following condition holds: 

\begin{align} \label{eq:weak_admis}
    \text{adm}_{*} \bkt{X, Y} = \text{true} \iff \text{dim} \bkt{C_i \cap C_j} = 0,
\end{align}
i.e., the intersection of $C_i$ and $C_j$ is at most a vertex.

The interaction list of a cluster $X$ is defined as
\begin{align*} \label{eq:IL_weak}
    \mathcal{IL}_{*} \bkt{X} := \big \{Y : \text{adm}_{*} \bkt{X, Y} = \text{true} \text{ and } \text{adm}_{*} \bkt{\text{parent} \bkt{X}, \text{parent} \bkt{Y}} = \text{false} \big \}.
\end{align*}

The neighbors list of $X$ is given by
\begin{align*} \label{eq:N_weak}
    \mathcal{N}_{*} \bkt{X} := \big \{ Y : \text{adm}_{*} \bkt{X, Y} = \text{false} \text{ and } Y \neq X \big \}.
\end{align*}

\begin{figure}[H]
\begin{center}
\subfloat[Level $1$]{
    \begin{tikzpicture}[scale=0.8]
        \fill [red] (1,0) rectangle (2,1);
        \fill [orange] (0,0) rectangle (1,1);
        \fill [orange] (1,1) rectangle (2,2);
        \fill [cyan] (0,1) rectangle (1,2);
        \draw[step=1cm, black] (0, 0) grid (2, 2);
        % \node[anchor=north] at (1.5,0.7) {\tiny $\mathcal{C}^{(1)}$};
        \draw[line width=0.5mm,  black] (0, 0) rectangle (2, 2);
    \end{tikzpicture}
    \label{hodlr2d_interaction_1a}
    }\quad
    \subfloat[Level $2$]{
    \begin{tikzpicture}[scale=0.4]
        \draw[black] (0, 0) grid (4, 4);
        \fill [green] (0,0) rectangle (1,1);
        \fill [green] (0,1) rectangle (1,2);
        \fill [green] (2,3) rectangle (3,4);
        \fill [green] (3,3) rectangle (4,4);
        \fill [red] (2,1) rectangle (3,2);
        \fill [cyan] (1,0) rectangle (2,1);
        \fill [orange] (1,1) rectangle (2,2);
        % \fill [cyan] (1,2) rectangle (2,3);
        \fill [orange] (2,2) rectangle (3,3);
        \fill [cyan] (3,2) rectangle (4,3);
        \fill [orange] (2,0) rectangle (3,1);
        \fill [cyan] (3,0) rectangle (4,1);
        \fill [orange] (3,1) rectangle (4,2);
        % \node[anchor=north] at (2.5,1.8) {\tiny $\mathcal{C}^{(2)}$};
        \node[anchor=north] at (1.5,2.8) {};
        \node[anchor=north] at (1.5,1.8) {};
        \node[anchor=north] at (.5,1.8) {};
        \draw[black] (0, 0) grid (4, 4);
        \fill [white] (0,2) rectangle (2,4);
        \draw[line width=0.5mm,  black] (0, 0) rectangle (2, 2);
        \draw[line width=0.5mm,  black] (2, 2) rectangle (4, 4);
        \draw[line width=0.5mm,  black] (2, 0) rectangle (4, 2);
    \end{tikzpicture}
    \label{hodlr2d_interaction_2a}
    }\quad  
    \subfloat[Level $3$]{
    \begin{tikzpicture}[scale=0.4]
        \draw[black] (0, 0) grid (4, 4);
        % \fill [green!50] (0,0) rectangle (1,1);
        % \fill [green!50] (0,1) rectangle (1,2);
        % \fill [green!50] (2,3) rectangle (3,4);
        % \fill [green!50] (3,3) rectangle (4,4);
        % \fill [cyan!50] (1,0) rectangle (2,1);
        % % \fill [cyan] (1,2) rectangle (2,3);
        % \fill [cyan!50] (3,2) rectangle (4,3);
        % \fill [cyan!50] (3,0) rectangle (4,1);
        \fill [red] (2,1) rectangle (3,2);
        \fill [green] (1,1) rectangle (2,2);
        \fill [orange] (1.5,1.5) rectangle (2,2);
        \fill [orange] (2.5,1.5) rectangle (3,2);
        \fill [cyan] (2.5,1) rectangle (3,1.5);
        \fill [orange] (2,1) rectangle (2.5,1.5);
        \fill [orange] (2,2) rectangle (2.5,2.5);
        \fill [green] (2,2) rectangle (3,3);
        \fill [green] (2,0) rectangle (3,1);
        \fill [green] (3,1) rectangle (4,2);
        % \node[anchor=north] at (2.28,2.1) {\tiny $\mathcal{C}^{(3)}$};
        \node[anchor=north] at (1.5,2.8) {};
        \node[anchor=north] at (1.5,1.8) {};
        \node[anchor=north] at (.5,1.8) {};
        \draw[black] (0, 0) grid (4, 4);
        \draw[step=0.5cm, black] (3, 1) grid (4, 2);
        \draw[step=0.5cm, black] (2, 0) grid (3, 1);
        \draw[step=0.5cm, black] (2, 2) grid (3, 3);
        \draw[step=0.5cm, black] (2, 1) grid (3, 2);
        \draw[step=0.5cm, black] (1, 1) grid (2, 2);
        \fill [white] (0,2) rectangle (2,4);
        \fill [orange] (2,2) rectangle (2.5,2.5);
        \fill [cyan] (2.5,2) rectangle (3,2.5);
        \fill [cyan] (1.5,1) rectangle (2,1.5);
        \draw[line width=0.5mm,  black] (1, 1) rectangle (2, 2);
        \draw[line width=0.5mm,  black] (2, 2) rectangle (3, 3);
        \draw[line width=0.5mm,  black] (3, 1) rectangle (4, 2);
        \draw[line width=0.5mm,  black] (2, 0) rectangle (3, 1);
        \draw[black] (0, 2) rectangle (2, 4);
        \draw[step=1cm, black] (0, 0) grid (2, 2);
    \end{tikzpicture}
    \label{hodlr2d_interaction_3a}
    }\quad
    \subfloat{
        \begin{tikzpicture}
            [%%%%%%%%%%%%%%%%%%%%%%%%%%%%%%
            box/.style={rectangle,draw=black, minimum size=0.2cm},scale=0.1
            ]%%%%%%%%%%%%%%%%%%%%%%%%%%%%%%
            \node[box,fill=red,,font=\tiny,label=right:Cluster considered,  anchor=west] at (8,0){};
            \node[box,fill=orange,,font=\tiny,label=right:Neighbor (edge-sharing) cluster,  anchor=west] at (8,4){};
            \node[box,fill=cyan,,font=\tiny,label=right:Vertex-sharing cluster,  anchor=west] at (8,8){};
            \node[box,fill=green,,font=\tiny,label=right:Far-field cluster,  anchor=west] at (8,12){};      
        \end{tikzpicture}
    }
    \caption{The interaction and neighbor lists of a cluster under our weak admissibility condition at different levels in $2$D. The admissible clusters are the green and cyan colored clusters enclosed within a noticeable \textbf{black} border.}
    \label{fig:HODLRdD_interaction_a}
\end{center}
\end{figure}
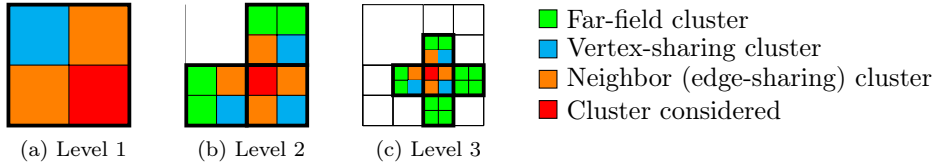

The interaction list and the neighbors list under the weak admissibility condition \cref{eq:weak_admis} are shown in \Cref{fig:HODLRdD_interaction_a}.

An upper bound on the size of the interaction list and the neighbor list of a cluster under weak admissibility condition \cref{eq:weak_admis} can be derived as follows \cite{nHODLRdD}
\begin{align} \label{eq:HODLRdD_matrix_IL_N}
   \abs{\mathcal{IL}_{*} \bkt{X}} \leq \bkt{2^d -1} \bkt{3^d - 2^d}, \qquad \abs{\mathcal{N}_{*} \bkt{X}} \leq \bkt{3^d - 2^d -1}.
\end{align}

\begin{remark}
The weak admissibility condition \cref{eq:weak_admis} can not be obtained from the standard admissibility condition \cref{eq:H_admiss} as the distance between two vertex-sharing clusters is zero (they share a common point). More precisely, let $X$ and $Y$ be two clusters at the same level of the tree with $Y \in \mathcal{IL}_{*} \bkt{X}$. If the corresponding hypercubes $C_i$ and $C_j$ share a vertex, then $\text{dist}\bkt{C_i,C_j} = 0$. Since $\text{dist}\bkt{C_i,C_j} = 0$, no real value of $\eta$ can make the weak admissibility condition \cref{eq:weak_admis} follow from \cref{eq:H_admiss}.
\end{remark}

\subsection{\texorpdfstring{$\mathcal{H}$-matrices}{H-matrices}}
Hierarchical matrices are multilevel block low-rank representations in which admissible blocks are compressed. More precisely, if $X \times Y$ is an admissible block, then the block matrix $K_{X,Y} \in \mathbb{C}^{m \times n}$ is approximated in low-rank form as follows:
\begin{equation*}
     K_{X,Y} \approx U V^*, \qquad \text{ where } \abs{X} = m, \space \abs{Y} = n \text{ and } U, V \in \mathbb{C}^{m \times r}.
\end{equation*} 
We rely on algebraic low-rank approximation techniques in this work, such as adaptive cross approximation (ACA) \cite{aca,bebendorf2012constructing}. In the $\mathcal{H}$-matrices framework, the bases associated with admissible blocks are constructed independently (in a non-nested manner).

\begin{figure}[H]
    \centering
    \subfloat[Level $1$]{
        \includegraphics[scale=0.38]{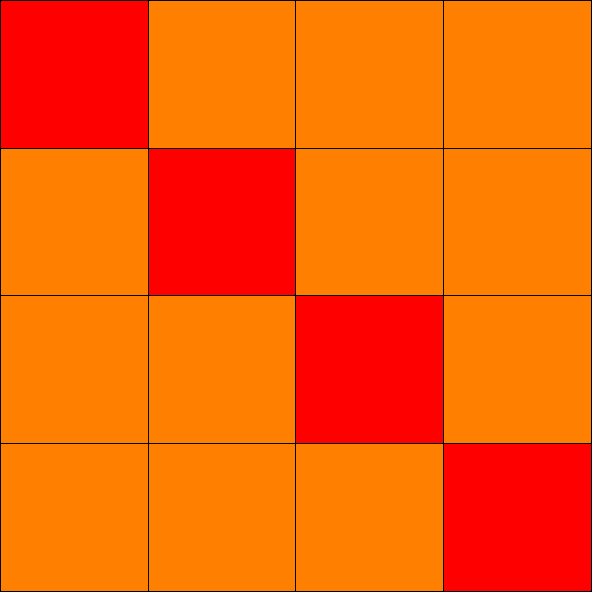}
        \label{fmm_mat_1}
        }
    \subfloat[Level $2$]{
        \includegraphics[scale=0.38]{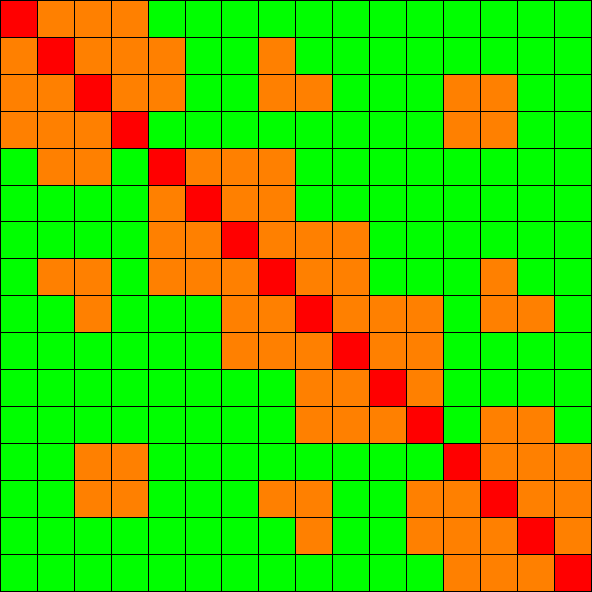}
        \label{fmm_mat_2}
        }
    \subfloat[Level $3$]{
        \includegraphics[scale=0.38]{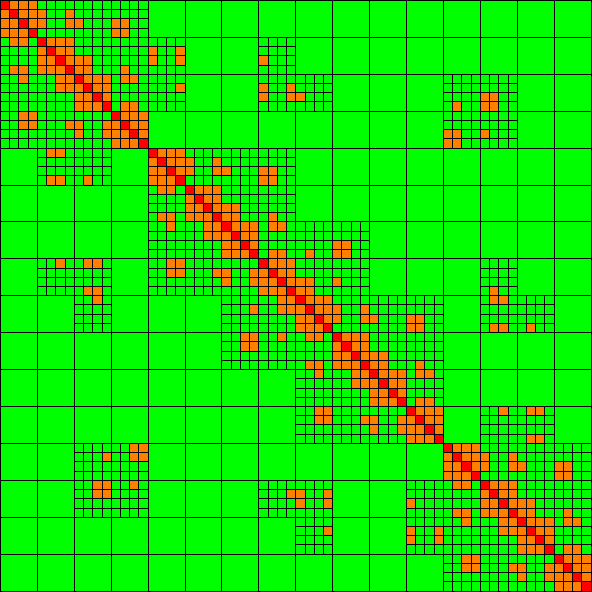}
        \label{fmm_mat_3}
        } \quad
    \subfloat{
        \begin{tikzpicture}
            [%%%%%%%%%%%%%%%%%%%%%%%%%%%%%%
            box/.style={rectangle,draw=black, minimum size=0.1cm},scale=0.1
            ]%%%%%%%%%%%%%%%%%%%%%%%%%%%%%%
            \node[box,fill=red,,font=\tiny,label=right:Diagonal block (dense),  anchor=west] at (-20,8){};
            \node[box,fill=orange,,font=\tiny,label=right:Neighbor block (dense),  anchor=west] at (21,8){};
            \node[box,fill=green,,font=\tiny,label=right:Far-field compressed block,  anchor=west] at (62,8){};
        \end{tikzpicture}
    }
    \caption{Hierarchical matrix based on standard admissibility condition \cref{eq:H_admiss} with $\eta = \sqrt{d}$ at different levels in $2$D (no compression at the level $1$).}
    \label{fig:FMM2D_matrix}
\end{figure}

\begin{figure}
    \centering
    \subfloat[Level $1$]{
        \includegraphics[scale=0.38]{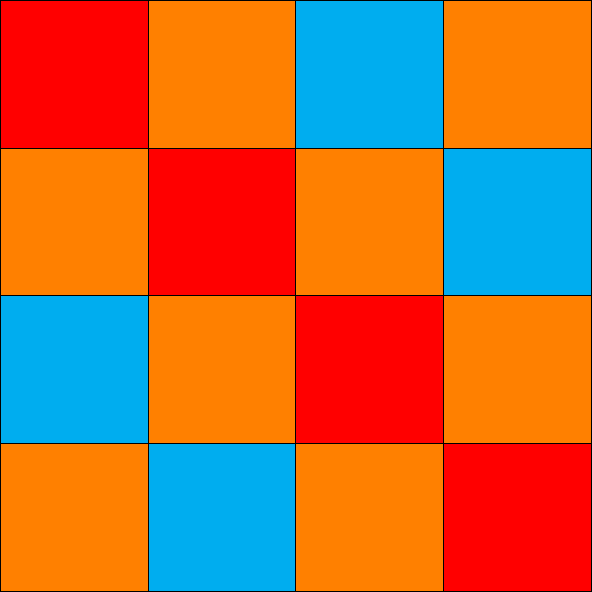}
        \label{fig:hodlr2d_mat_1a}
        }
    \subfloat[Level $2$]{
        \includegraphics[scale=0.38]{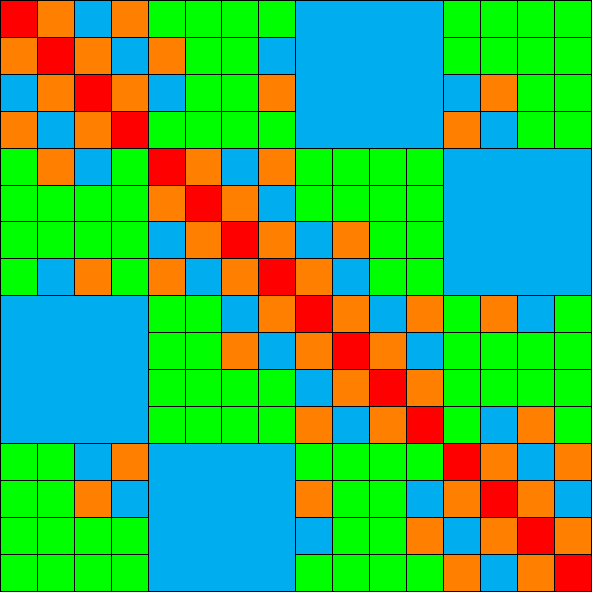}
        \label{fig:hodlr2d_mat_2a}
        }
    \subfloat[Level $3$]{
        \includegraphics[scale=0.38]{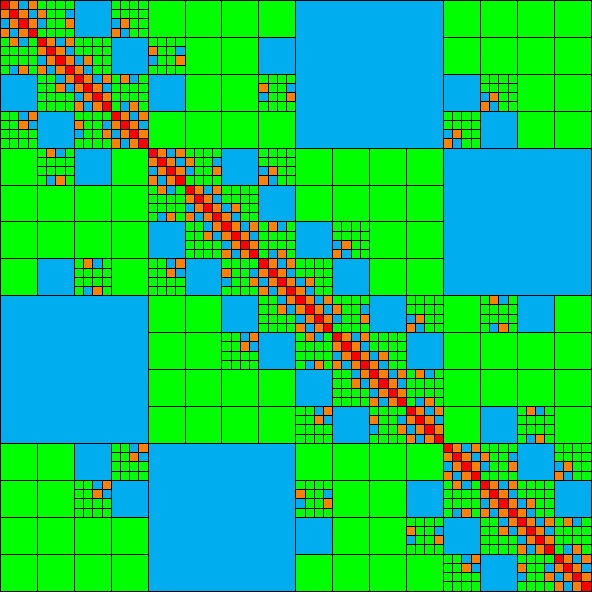}
        \label{fig:hodlr2d_mat_3a}
        }\quad 
    \subfloat{
        \begin{tikzpicture}
            [%%%%%%%%%%%%%%%%%%%%%%%%%%%%%%
            box/.style={rectangle,draw=black, minimum size=0.1cm},scale=0.1
            ]%%%%%%%%%%%%%%%%%%%%%%%%%%%%%%
            \node[box,fill=red,,font=\tiny,label=right:Diagonal block (dense),  anchor=west] at (0,16){};
            \node[box,fill=orange,,font=\tiny,label=right:Neighbor block (dense),  anchor=west] at (0,12){};
            \node[box,fill=cyan,,font=\tiny,label=right:Vertex-sharing compressed block,  anchor=west] at (50,16){};
            \node[box,fill=green,,font=\tiny,label=right:Far-field compressed block,  anchor=west] at (50,12){};
        \end{tikzpicture}
    }
    \caption{Hierarchical matrix based on our weak admissibility condition \cref{eq:weak_admis} at different levels in $2$D (no far-field compression at the level $1$).}
    \label{fig:HODLRdD_matrix_a}
\end{figure}
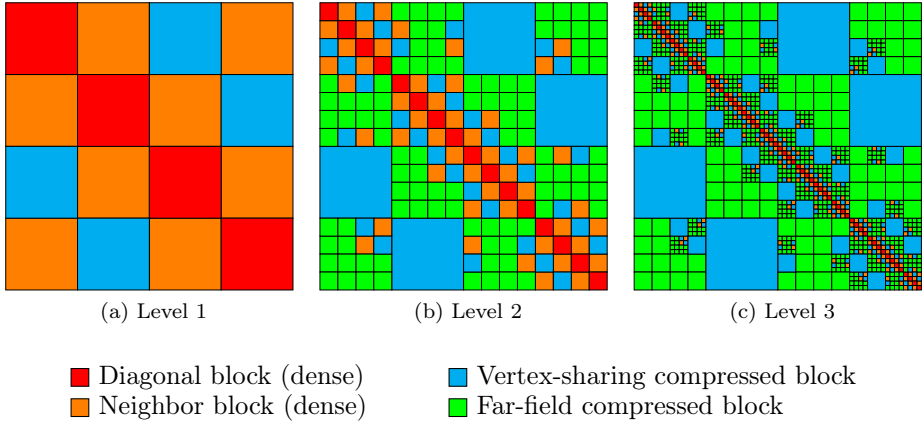

The hierarchical matrices under the standard \cref{eq:H_admiss} and weak admissibility \cref{eq:weak_admis} conditions are illustrated in \Cref{fig:FMM2D_matrix} and \Cref{fig:HODLRdD_matrix_a}, respectively. A level-wise comparison helps in understanding the resulting block structures.

\subsubsection{Adaptive cross approximation}
The admissible block matrix $K_{X,Y}$ can be approximated by ACA \cite{aca,bebendorf2012constructing} as follows:
\begin{equation}
   K_{X, Y} \approx K_{X, Y}^{(r)} =  UV^* = K_{X,\sigma_Y} \bkt{K_{\tau_X,\sigma_Y}}^{-1} K_{\tau_X,Y}
\end{equation}
where $\tau_X \subset X$ and $\sigma_Y \subset Y$ are called the \emph{row} and \emph{column} pivots, respectively. The matrix $K_{X, Y}^{(r)}$ is the $r^{th}$ update of the original matrix $K_{X, Y}$. Here, we employ the partially pivoted ACA algorithm \cite{aca}. The partially pivoted ACA algorithm stops its iterations for a prescribed tolerance $\epsilon$ when the following condition is satisfied
\begin{equation}
    \magn{u_r}_2 \magn{v_r}_2 \leq \epsilon \magn{K_{X, Y}^{(r)}}_F
\end{equation}
where $u_r$ and $v_r$ are the $r^{th}$ column vectors of the matrices $U$ and $V$, respectively.

\subsection{\texorpdfstring{$\mathcal{H}^2$-matrices}{H2-matrices}}
The $\mathcal{H}$-matrix framework can be further improved by employing nested cluster bases, resulting in the $\mathcal{H}^2$-matrix format. We construct the $\mathcal{H}^2$-matrix representation using nested cross approximation (NCA) \cite{bebendorf2012constructing, zhao2019fast}.

\subsubsection{Nested cross approximation} \label{sec:nca}
NCA can be viewed as an extension of ACA \cite{aca}, in which carefully selected matrix indices are used to generate representative bases while preserving their nested structure.

The low-rank approximation of an admissible block/sub-matrix via NCA takes the following form \cite{bebendorf2012constructing, zhao2019fast}

\begin{equation} \label{nca_equ}
    K_{X, Y} \approx \underbrace{K_{X, \sigma_X} \bkt{K_{\tau_X, \sigma_X}}^{-1}}_{\mathcal{U}_X} \underbrace{K_{\tau_{X}, \sigma_{Y}}}_{\mathcal{K}_{X,Y}} \underbrace{\bkt{K_{\tau_{Y}, \sigma_{Y}}}^{-1} K_{\tau_{Y}, {Y}}}_{\mathcal{V}_{Y}^*},
\end{equation}
where $\tau_{X} \subset X$, $\sigma_{X} \subset \mathcal{A} \bkt{X}$ and $\tau_{Y} \subset \mathcal{A} \bkt{Y}$, $\sigma_{Y} \subset Y$. The set $\mathcal{A} \bkt{\mathcal{C}}$ represent the union of all the clusters forming the interaction list at the same level and the clusters that form the interaction list with $\mathcal{C}$'s ancestors at the lower (coarser) level of the tree.

Although this work focuses on developing new hierarchical matrices, we adopt the FMM terminology for convenience, such as the translation operators (M$2$M, M$2$L, and L$2$L). 

The translation operators are obtained from the computed pivots as follows.

\begin{itemize}
    \item \textbf{L$2$P and P$2$M operators.} For a \emph{leaf level} admissible block $X \times Y$, the matrices $\mathcal{U}_X$ and $\mathcal{V}_{Y}^*$ are given by
\end{itemize}
    \begin{equation} \label{L2P_operator}
        \mathcal{U}_X = K_{X,\sigma_X} \bkt{K_{\tau_X, \sigma_X}}^{-1} 
    \end{equation}

    \begin{equation} \label{P2M_operator}
        \mathcal{V}_{Y}^* = \bkt{K_{\tau_Y, \sigma_Y}}^{-1} K_{\tau_Y, Y}
    \end{equation}
The matrices $\mathcal{U}_X$ and $\mathcal{V}_{Y}^*$ are called \textbf{L$2$P} (local to particles) and \textbf{P$2$M} (particles to multipole) operators, respectively.

\begin{itemize}
    \item \textbf{L$2$L and M$2$M operators.} For a \emph{non-leaf level} admissible block $X \times Y$, the matrices $\mathcal{U}_X$ and $\mathcal{V}_{Y}^*$ are given by 
\end{itemize}
\begin{equation} \label{L2L_operator}
                \mathcal{U}_{X} = \begin{bmatrix}
                    \mathcal{U}_{X_{1}} & 0 & \hdots & 0\\
                    0 & \mathcal{U}_{X_{2}} & \hdots & 0\\
                    \vdots & \vdots & \ddots & \vdots\\
                    0 & 0 & \hdots & \mathcal{U}_{X_{2^{d}}}
                    \end{bmatrix}
                    \begin{bmatrix}
                    \Tilde{\mathcal{U}}_{X_{1} X}\\
                    \Tilde{\mathcal{U}}_{X_{2} X}\\
                    \vdots\\
                    \Tilde{\mathcal{U}}_{X_{2^{d}} X}\\
                    \end{bmatrix}
\end{equation}
where $ X = \displaystyle \bigcup_{c=1}^{2^d} X_c$ with $X_c$ is a child of $X$ and
\begin{equation} 
    \Tilde{\mathcal{U}}_{X_{c} X} = K_{\tau_{X_c},\sigma_{X}} \bkt{K_{\tau_{X}, \sigma_X}}^{-1}, \qquad  1 \leq c \leq 2^d.
\end{equation}

The matrices $\Tilde{\mathcal{U}}_{X_{c} X}$ are called the \textbf{L$2$L} (local to local) operators.

\begin{equation} \label{M2M_operator}
                \mathcal{V}_{Y} = \begin{bmatrix}
                    \mathcal{V}_{Y_{1}} & 0 & \hdots & 0\\
                    0 & \mathcal{V}_{Y_{2}} & \hdots & 0\\
                    \vdots & \vdots & \ddots & \vdots\\
                    0 & 0 & \hdots & \mathcal{V}_{Y_{2^{d}}}
                    \end{bmatrix}
                    \begin{bmatrix}
                    \Tilde{\mathcal{V}}_{{Y Y_{1}}^*}\\
                    \Tilde{\mathcal{V}}_{{Y Y_{2}}^*}\\
                    \vdots\\
                    \Tilde{\mathcal{V}}_{{Y Y_{2^{d}}}^*}\\
                    \end{bmatrix}
\end{equation}
where $ Y = \displaystyle \bigcup_{c=1}^{2^d} Y_c$ with $Y_c$ is a child of $Y$ and
\begin{equation}
  \Tilde{\mathcal{V}}_{{Y Y_{c}}^*} = \bkt{K_{\tau_Y, \sigma_Y}}^{-1} K_{\tau_Y, \sigma_{Y_c}}, \qquad  1 \leq c \leq 2^d.
\end{equation}

The matrices $\Tilde{\mathcal{V}}_{{Y Y_{c}}^*}$ are called the \textbf{M$2$M} (multipole to multipole) operators. 

\begin{itemize}
    \item  \textbf{M$2$L operators.} For any admissible block $X \times Y$, the \textbf{M$2$L} (multipole to local) operator is given by
\end{itemize}
    \begin{equation}
        \mathcal{K}_{X,Y} = K_{\tau_X, \sigma_Y}, \quad Y \in \mathcal{IL} \bkt{X}.
    \end{equation}

\begin{remark}
    If the kernel matrix is symmetric, the P$2$M or M$2$M operators can be obtained simply by taking the transpose of L$2$P or L$2$L operators, respectively. 
\end{remark}

The validation of NCA, along with numerical error analysis, is presented in \cite{bebendorf2012constructing}. The efficiency of NCA depends on the strategy employed for selecting the pivots that represent the blocks at a specific level of the tree. 

Bebendorf et al. \cite{bebendorf2012constructing} choose indices of points close to the Chebyshev nodes and perform ACA upon them to obtain the pivots. The tree is traversed in a top to bottom manner, resulting in an overall construction complexity of $\mathcal{O} \bkt{N \log N}$. 

Zhao et al. \cite{zhao2019fast} eliminate the geometric dependence by introducing two purely algebraic variants of NCA. They first develop an $\mathcal{O}(N\log N)$ algorithm for pivot selection that employs a top to bottom tree traversal, analogous to the approach in \cite{bebendorf2012constructing}, but does not require geometric information about the clusters. In their second approach, local pivots are identified via a bottom to top tree traversal and subsequently refined through a top to bottom traversal to obtain global pivots. This procedure reduces the overall construction complexity to $\mathcal{O}\bkt{N}$. 

Recently, Gujjula et al. \cite{gujjula2022new} demonstrated that the global pivot selection step of \cite{zhao2019fast} could be omitted without affecting the approximation quality. Their NCA construction relies solely on the local pivots obtained during the bottom to top tree traversal and retains $\mathcal{O}\bkt{N}$ complexity. Their numerical experiments further show that this simplification can provide a computational advantage over the earlier NCA variants \cite{bebendorf2012constructing, zhao2019fast} without compromising accuracy.

\section{New hierarchical matrices} \label{sec:hybrid_matrix}
This section discusses two new hierarchical matrices to perform fast matrix-vector product: $\mathcal{H}_{hw}^{\prime}$, $\mathcal{H}_{hw}^{\prime \prime}$. The proposed hierarchical matrices are based on the following novel hybrid weak admissibility condition. 

\subsection*{Hybrid weak admissibility condition} \label{sec:hybrid_admiss}
One drawback of $\mathcal{H}$-matrices is the large number of admissible blocks as tree depth increases, along with many dense blocks at the leaf level. In multidimensional problems, the rapidly growing number of admissible blocks can become a computational bottleneck as the tree depth increases. To mitigate this drawback, while keeping the leaf-level block size (or the tree depth $L$) unchanged, we treat all the adjacent clusters as admissible at finer levels. This is motivated by the observation that admissible blocks at these levels typically have relatively low numerical rank due to their smaller block sizes.

Let $\mathcal{C}_i$ and $\mathcal{C}_j$ be a pair of clusters at level $l$ of the balanced $2^d$-tree, and let $\ell < L$ denote an intermediate level of the tree, which we referred to as the switching level. We define the \emph{hybrid weak admissibility condition} as follows:

\begin{align} \label{eq:hybrid_admiss}
       \text{adm}_{hw} \bkt{\mathcal{C}_i, \mathcal{C}_j} = \text{true} \iff
\begin{cases}
    \text{adm}_* \bkt{\mathcal{C}_i, \mathcal{C}_j} = \text{true} ,& \text{if } 1 \leq l \leq \ell \\
    \text{adm}_* \bkt{\mathcal{C}_i, \mathcal{C}_j} = \text{false} \text{ and } i \neq j,& \text{if } l = \ell < L \\
    \mathcal{C}_i \in siblings(\mathcal{C}_j) \text{ and } i \neq j ,& \text{if } \ell +1 \leq l \leq L.
\end{cases}
\end{align}

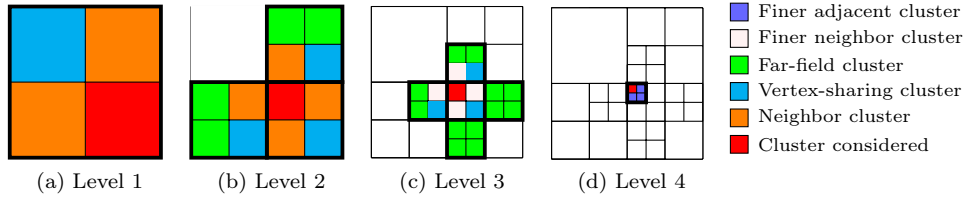
\begin{figure}
\begin{center}
\subfloat[Level $1$]{
    \begin{tikzpicture}[scale=1]
        \fill [red] (1,0) rectangle (2,1);
        \fill [orange] (0,0) rectangle (1,1);
        \fill [orange] (1,1) rectangle (2,2);
        \fill [cyan] (0,1) rectangle (1,2);
        \draw[step=1cm, black] (0, 0) grid (2, 2);
        % \node[anchor=north] at (1.5,0.7) {\tiny $\mathcal{C}^{(1)}$};
        \draw[line width=0.5mm,  black] (0, 0) rectangle (2, 2);
    \end{tikzpicture}
    \label{hodlr2d_interaction_1}
    }
    \subfloat[Level $2$]{
    \begin{tikzpicture}[scale=0.5]
        \draw[black] (0, 0) grid (4, 4);
        \fill [green] (0,0) rectangle (1,1);
        \fill [green] (0,1) rectangle (1,2);
        \fill [green] (2,3) rectangle (3,4);
        \fill [green] (3,3) rectangle (4,4);
        \fill [red] (2,1) rectangle (3,2);
        \fill [cyan] (1,0) rectangle (2,1);
        \fill [orange] (1,1) rectangle (2,2);
        % \fill [cyan] (1,2) rectangle (2,3);
        \fill [orange] (2,2) rectangle (3,3);
        \fill [cyan] (3,2) rectangle (4,3);
        \fill [orange] (2,0) rectangle (3,1);
        \fill [cyan] (3,0) rectangle (4,1);
        \fill [orange] (3,1) rectangle (4,2);
        % \node[anchor=north] at (2.5,1.8) {\tiny $\mathcal{C}^{(2)}$};
        \node[anchor=north] at (1.5,2.8) {};
        \node[anchor=north] at (1.5,1.8) {};
        \node[anchor=north] at (.5,1.8) {};
        \draw[black] (0, 0) grid (4, 4);
        \fill [white] (0,2) rectangle (2,4);
        \draw[line width=0.5mm,  black] (0, 0) rectangle (2, 2);
        \draw[line width=0.5mm,  black] (2, 2) rectangle (4, 4);
        \draw[line width=0.5mm,  black] (2, 0) rectangle (4, 2);
    \end{tikzpicture}
    \label{hodlr2d_interaction_2}
    }
    \subfloat[Level $3$]{
    \begin{tikzpicture}[scale=0.5]
        \draw[black] (0, 0) grid (4, 4);
        \draw[black] (1, 0) grid (4, 4);
        \fill [white] (0,0) rectangle (1,1);
        \fill [white] (0,1) rectangle (1,2);
        \fill [white] (2,3) rectangle (3,4);
        \fill [white] (3,3) rectangle (4,4);
        \fill [white] (1,0) rectangle (2,1);
        \fill [white] (3,2) rectangle (4,3);
        \fill [white] (3,0) rectangle (4,1);
        \fill [red] (2,1) rectangle (3,2);
        \fill [green] (1,1) rectangle (2,2);
        \fill [pink!20] (1.5,1.5) rectangle (2,2);
        \fill [pink!20] (2.5,1.5) rectangle (3,2);
        \fill [cyan] (2.5,1) rectangle (3,1.5);
        \fill [pink!20] (2,1) rectangle (2.5,1.5);
        \fill [pink!20] (2,2) rectangle (2.5,2.5);
        \fill [green] (2,2) rectangle (3,3);
        \fill [green] (2,0) rectangle (3,1);
        \fill [green] (3,1) rectangle (4,2);
        % \node[anchor=north] at (2.28,2.1) {\tiny $\mathcal{C}^{(3)}$};
        \node[anchor=north] at (1.5,2.8) {};
        \node[anchor=north] at (1.5,1.8) {};
        \node[anchor=north] at (.5,1.8) {};
        \draw[black] (0, 0) grid (4, 4);
        \draw[step=0.5cm, black] (3, 1) grid (4, 2);
        \draw[step=0.5cm, black] (2, 0) grid (3, 1);
        \draw[step=0.5cm, black] (2, 2) grid (3, 3);
        \draw[step=0.5cm, black] (2, 1) grid (3, 2);
        \draw[step=0.5cm, black] (1, 1) grid (2, 2);
        \fill [white] (0,2) rectangle (2,4);
        \fill [pink!20] (2,2) rectangle (2.5,2.5);
        \fill [cyan] (2.5,2) rectangle (3,2.5);
        \fill [cyan] (1.5,1) rectangle (2,1.5);
        \draw[line width=0.5mm,  black] (1, 1) rectangle (2, 2);
        \draw[line width=0.5mm,  black] (2, 2) rectangle (3, 3);
        \draw[line width=0.5mm,  black] (3, 1) rectangle (4, 2);
        \draw[line width=0.5mm,  black] (2, 0) rectangle (3, 1);
        \draw[black] (0, 2) rectangle (2, 4);
        \draw[step=1cm, black] (0, 0) grid (2, 2);
    \end{tikzpicture}
    \label{hodlr2d_interaction_3}
    }
    \subfloat[Level $4$]{
    \begin{tikzpicture}[scale=0.5]
        \draw[black] (0, 0) grid (4, 4);
        \draw[black] (1, 0) grid (4, 4);
        \fill [white] (0,0) rectangle (1,1);
        \fill [white] (0,1) rectangle (1,2);
        \fill [white] (2,3) rectangle (3,4);
        \fill [white] (3,3) rectangle (4,4);
        \fill [white] (1,0) rectangle (2,1);
        \fill [white] (3,2) rectangle (4,3);
        \fill [white] (3,0) rectangle (4,1);
        \fill [red] (2,1) rectangle (3,2);
        \fill [white] (1,1) rectangle (2,2);
        \fill [white] (1.5,1.5) rectangle (2,2);
        \fill [white] (2.5,1.5) rectangle (3,2);
        \fill [white] (2.5,1) rectangle (3,1.5);
        \fill [white] (2,1) rectangle (2.5,1.5);
        \fill [white] (2,2) rectangle (2.5,2.5);
        \fill [white] (2,2) rectangle (3,3);
        \fill [white] (2,0) rectangle (3,1);
        \fill [white] (3,1) rectangle (4,2);
        % \node[anchor=north] at (2.28,2.1) {\tiny $\mathcal{C}^{(3)}$};
        \node[anchor=north] at (1.5,2.8) {};
        \node[anchor=north] at (1.5,1.8) {};
        \node[anchor=north] at (.5,1.8) {};
        \draw[black] (0, 0) grid (4, 4);
        \fill [white] (0,2) rectangle (2,4);
        \draw[step=0.5cm, black] (3, 1) grid (4, 2);
        \draw[step=0.5cm, black] (2, 0) grid (3, 1);
        \draw[step=0.5cm, black] (2, 2) grid (3, 3);
        \draw[step=0.5cm, black] (2, 1) grid (3, 2);
        \draw[step=0.5cm, black] (1, 1) grid (2, 2);
        % \fill [white] (0,2) rectangle (2,4);
        % \fill [white] (2,2) rectangle (2.5,2.5);
        % \fill [white] (2.5,2) rectangle (3,2.5);
        \draw[black] (2, 2) grid (2.5, 2.5);
        % \fill [cyan!20] (1.5,1) rectangle (2,1.5);
        \fill [blue!60] (2.25,1.75) rectangle (2.5,2);
        \fill [blue!60] (2,1.5) rectangle (2.25,1.75);
        \fill [blue!60] (2.25,1.5) rectangle (2.5,1.75);
        \draw[line width=0.5mm,  black] (2, 1.5) rectangle (2.5, 2);
        \draw[step=0.25cm, black] (2, 1.5) grid (2.5, 2);
        \draw[black] (0, 2) rectangle (2, 4);
        \draw[step=1cm, black] (0, 0) grid (2, 2);
    \end{tikzpicture}
    \label{hodlr2d_interaction_4}
    }
    \subfloat{
        \begin{tikzpicture}
            [%%%%%%%%%%%%%%%%%%%%%%%%%%%%%%
            box/.style={rectangle,draw=black, minimum size=0.2cm},scale=0.1
            ]%%%%%%%%%%%%%%%%%%%%%%%%%%%%%%
            \node[box,fill=red,,font=\tiny,label=right:\scriptsize Cluster considered,  anchor=west] at (8,-0.5){};
            \node[box,fill=orange,,font=\tiny,label=right:\scriptsize Neighbor cluster,  anchor=west] at (8,3){};
            \node[box,fill=cyan,,font=\tiny,label=right:\scriptsize Vertex-sharing cluster,  anchor=west] at (8,6.5){};
            \node[box,fill=green,,font=\tiny,label=right:\scriptsize Far-field cluster,  anchor=west] at (8,10){};  
            \node[box,fill=pink!20,,font=\tiny,label=right:\scriptsize Finer neighbor cluster,  anchor=west] at (8,13.5){};  
            \node[box,fill=blue!60,,font=\tiny,label=right:\scriptsize Finer adjacent cluster,  anchor=west] at (8,17){};  
        \end{tikzpicture}
    }
    \caption{Interaction list under the hybrid weak admissibility condition \cref{eq:hybrid_admiss} with $\ell=3$ and $L=4$. At level $3$, finer neighbor clusters (pink color) are considered admissible. At level $4$, adjacent clusters (blue color) are considered admissible.}
    \label{fig:HODLRdD_interaction}
\end{center}
\end{figure}

\begin{figure}
    \centering
    \subfloat[\scriptsize Level $1$]{
        \includegraphics[scale=0.55]{Matrix_images/HODLR2D_1.pdf}
        \label{mixed_hodlr2d_mat_1}
        }\qquad
    \subfloat[\scriptsize Level $2$]{
        \includegraphics[scale=0.55]{Matrix_images/HODLR2D_2.pdf}
        \label{mixed_hodlr2d_mat_2}
        }\qquad 
    \subfloat[\scriptsize Level $3$]{
        \includegraphics[scale=0.55]{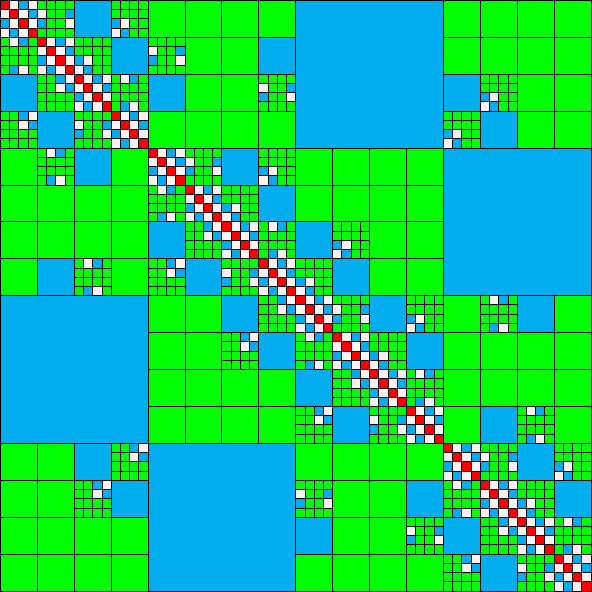}
        \label{mixed_hodlr2d_mat_3}
        }\qquad 
    \subfloat[\scriptsize Level $4$]{
        \includegraphics[scale=0.55]{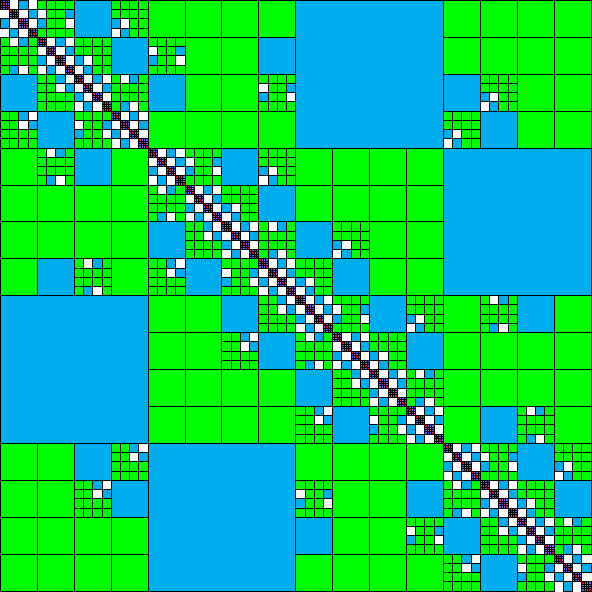}
        \label{mixed_hodlr2d_mat_4}
        }\qquad 
    \subfloat{
        \begin{tikzpicture}
            [%%%%%%%%%%%%%%%%%%%%%%%%%%%%%%
            box/.style={rectangle,draw=black, minimum size=0.1cm},scale=0.1
            ]%%%%%%%%%%%%%%%%%%%%%%%%%%%%%%
            \node[box,fill=red,,font=\tiny,label=right:Dense diagonal block,  anchor=west] at (0,8){};
            \node[box,fill=orange,,font=\tiny,label=right:Dense neighbor block,  anchor=west] at (0,4){};
            \node[box,fill=pink,,font=\tiny,label=right:Compressed finer neighbor block,  anchor=west] at (0,0){};
            \node[box,fill=cyan,,font=\tiny,label=right:Compressed vertex-sharing block,  anchor=west] at (60,8){};
            \node[box,fill=green,,font=\tiny,label=right:Compressed far-field block,  anchor=west] at (60,4){};
            \node[box,fill=blue!60,,font=\tiny,label=right:Compressed finer adjacent block,  anchor=west] at (60,0){};
        \end{tikzpicture}
    }
    \caption{Hierarchical matrix under the hybrid weak admissibility condition at different levels. A \textbf{zoomed-in view} shows how the hierarchical matrices are constructed when $\ell=3$ and $L=4$.}
    \label{fig:mixed_HODLRdD_matrix}
\end{figure}

The hybrid admissibility condition employs the weak admissibility condition \cref{eq:weak_admis} at \emph{coarser levels}, where the block sizes are relatively \textbf{large}. At \emph{finer levels}, it treats all adjacent sibling clusters as admissible, where the block sizes are comparatively \textbf{small}. This approach reduces the number of admissible clusters at the \emph{finer levels}, and self-clusters are the only inadmissible clusters at the leaf level. The interaction list under the hybrid admissibility condition is illustrated in \Cref{fig:HODLRdD_interaction}. The corresponding hierarchical matrix is shown in \Cref{fig:mixed_HODLRdD_matrix}, where all blocks except the diagonal blocks are represented in compressed form.

We propose new hierarchical matrices under the hybrid weak admissibility condition \cref{eq:hybrid_admiss}, based on different bases associated with the admissible blocks. At \emph{coarser levels}, we employ both \textbf{nested} and \textbf{non-nested} bases, whereas at \emph{finer levels}, where adjacent clusters are admissible, we use only \textbf{non-nested} bases.

\subsection{\texorpdfstring{$\mathcal{H}_{hw}^{\prime}$}{H1}-matrix} \label{subsec:H1}
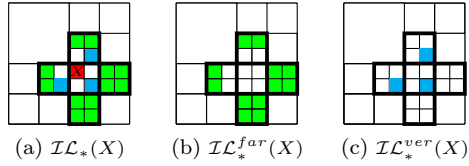
\begin{figure}[H]
\begin{center}
    \subfloat[\scriptsize $\mathcal{IL}_{*}(X)$]{
    \begin{tikzpicture}[scale=0.4]
        \draw[black] (0, 0) grid (4, 4);
        \fill [white] (0,0) rectangle (1,1);
        \fill [white] (0,1) rectangle (1,2);
        \fill [white] (2,3) rectangle (3,4);
        \fill [white] (3,3) rectangle (4,4);
        \fill [white] (1,0) rectangle (2,1);
        % \fill [cyan] (1,2) rectangle (2,3);
        \fill [white] (3,2) rectangle (4,3);
        \fill [white] (3,0) rectangle (4,1);
        \fill [red] (2,1) rectangle (3,2);
        \fill [green] (1,1) rectangle (2,2);
        \fill [white] (1.5,1.5) rectangle (2,2);
        \fill [white] (2.5,1.5) rectangle (3,2);
        \fill [cyan] (2.5,1) rectangle (3,1.5);
        \fill [white] (2,1) rectangle (2.5,1.5);
        \fill [white] (2,2) rectangle (2.5,2.5);
        \fill [green] (2,2) rectangle (3,3);
        \fill [green] (2,0) rectangle (3,1);
        \fill [green] (3,1) rectangle (4,2);
        \node[anchor=north] at (2.27,2.2) {\tiny $X$};
        \node[anchor=north] at (1.5,2.8) {};
        \node[anchor=north] at (1.5,1.8) {};
        \node[anchor=north] at (.5,1.8) {};
        \draw[black] (0, 0) grid (4, 4);
        \draw[step=0.5cm, black] (3, 1) grid (4, 2);
        \draw[step=0.5cm, black] (2, 0) grid (3, 1);
        \draw[step=0.5cm, black] (2, 2) grid (3, 3);
        \draw[step=0.5cm, black] (2, 1) grid (3, 2);
        \draw[step=0.5cm, black] (1, 1) grid (2, 2);
        \fill [white] (0,2) rectangle (2,4);
        \fill [white] (2,2) rectangle (2.5,2.5);
        \fill [cyan] (2.5,2) rectangle (3,2.5);
        \fill [cyan] (1.5,1) rectangle (2,1.5);
        \draw[line width=0.1mm,  black] (2.5, 2) rectangle (3, 2.5);
        \draw[line width=0.1mm,  black] (2, 2) rectangle (2.5, 2.5);
        \draw[line width=0.5mm,  black] (1, 1) rectangle (2, 2);
        \draw[line width=0.5mm,  black] (2, 2) rectangle (3, 3);
        \draw[line width=0.5mm,  black] (3, 1) rectangle (4, 2);
        \draw[line width=0.5mm,  black] (2, 0) rectangle (3, 1);
        \draw[black] (0, 2) rectangle (2, 4);
        \draw[step=1cm, black] (0, 0) grid (2, 2);
    \end{tikzpicture}
    }\quad
    \subfloat[\scriptsize $\mathcal{IL}_*^{far}(X)$]{
    \begin{tikzpicture}[scale=0.4]
        \draw[black] (0, 0) grid (4, 4);
        \fill [white] (0,0) rectangle (1,1);
        \fill [white] (0,1) rectangle (1,2);
        \fill [white] (2,3) rectangle (3,4);
        \fill [white] (3,3) rectangle (4,4);
        % \fill [cyan!25] (1,0) rectangle (2,1);
        % \fill [cyan] (1,2) rectangle (2,3);
        % \fill [cyan!25] (3,2) rectangle (4,3);
        % \fill [cyan!25] (3,0) rectangle (4,1);
        % \fill [red] (2,1) rectangle (3,2);
        \fill [green] (1,1) rectangle (2,2);
        \fill [white] (1.5,1.5) rectangle (2,2);
        \fill [white] (2.5,1.5) rectangle (3,2);
        \fill [white] (2.5,1) rectangle (3,1.5);
        \fill [white] (2,1) rectangle (2.5,1.5);
        \fill [white] (2,2) rectangle (2.5,2.5);
        \fill [green] (2,2) rectangle (3,3);
        \fill [green] (2,0) rectangle (3,1);
        \fill [green] (3,1) rectangle (4,2);
        \draw[step=0.5cm, black] (2, 0) grid (3, 1);
        % \node[anchor=north] at (2.27,2) {\tiny $\mathcal{C}$};
        \node[anchor=north] at (1.5,2.8) {};
        \node[anchor=north] at (1.5,1.8) {};
        \node[anchor=north] at (.5,1.8) {};
        \draw[black] (0, 0) grid (4, 4);
        \draw[step=0.5cm, black] (3, 1) grid (4, 2);
        \draw[step=0.5cm, black] (2, 0) grid (3, 1);
        \draw[step=0.5cm, black] (2, 2) grid (3, 3);
        \draw[step=0.5cm, black] (2, 1) grid (3, 2);
        \draw[step=0.5cm, black] (1, 1) grid (2, 2);
        % \fill [cyan!25] (0,2) rectangle (2,4);
        \fill [white] (2,2) rectangle (2.5,2.5);
        \fill [white] (2.5,2) rectangle (3,2.5);
        \draw[line width=0.1mm,  black] (2.5, 2) rectangle (3, 2.5);
        \draw[line width=0.1mm,  black] (2, 2) rectangle (2.5, 2.5);
        \fill [white] (1.5,1) rectangle (2,1.5);
        \fill [white] (0,2) rectangle (2,4);
        \draw[line width=0.5mm,  black] (1, 1) rectangle (2, 2);
        \draw[line width=0.5mm,  black] (2, 2) rectangle (3, 3);
        \draw[line width=0.5mm,  black] (3, 1) rectangle (4, 2);
        \draw[line width=0.5mm,  black] (2, 0) rectangle (3, 1);
        \draw[black] (0, 2) rectangle (2, 4);
        \draw[step=1cm, black] (0, 0) grid (2, 2);
        \draw[step=0.5cm, black] (1, 1) grid (2, 2);
    \end{tikzpicture}
    }\quad
    \subfloat[\scriptsize $\mathcal{IL}_*^{ver}(X)$]{
    \begin{tikzpicture}[scale=0.4]
        \draw[black] (0, 0) grid (4, 4);
        % \fill [green!50] (0,0) rectangle (1,1);
        % \fill [green!50] (0,1) rectangle (1,2);
        % \fill [green!50] (2,3) rectangle (3,4);
        % \fill [green!50] (3,3) rectangle (4,4);
        \fill [white] (1,0) rectangle (2,1);
        % \fill [cyan] (1,2) rectangle (2,3);
        \fill [white] (3,2) rectangle (4,3);
        \fill [white] (3,0) rectangle (4,1);
        % \fill [red] (2,1) rectangle (3,2);
        % \fill [green] (1,1) rectangle (2,2);
        % \fill [orange] (1.5,1.5) rectangle (2,2);
        % \fill [orange] (2.5,1.5) rectangle (3,2);
        \fill [cyan] (2.5,1) rectangle (3,1.5);
        % \fill [orange] (2,1) rectangle (2.5,1.5);
        % \fill [orange] (2,2) rectangle (2.5,2.5);
        % \fill [green] (2,2) rectangle (3,3);
        % \fill [green] (2,0) rectangle (3,1);
        % \fill [green] (3,1) rectangle (4,2);
        \node[anchor=north] at (1.5,2.8) {};
        \node[anchor=north] at (1.5,1.8) {};
        \node[anchor=north] at (.5,1.8) {};
        \draw[black] (0, 0) grid (4, 4);
        \draw[step=0.5cm, black] (3, 1) grid (4, 2);
        \draw[step=0.5cm, black] (2, 0) grid (3, 1);
        \draw[step=0.5cm, black] (2, 2) grid (3, 3);
        \draw[step=0.5cm, black] (2, 1) grid (3, 2);
        \draw[step=0.5cm, black] (1, 1) grid (2, 2);
        \fill [white] (0,2) rectangle (2,4);
        % \fill [orange] (2,2) rectangle (2.5,2.5);
        \fill [cyan] (2.5,2) rectangle (3,2.5);
        \fill [cyan] (1.5,1) rectangle (2,1.5);
        \draw[line width=0.5mm,  black] (1, 1) rectangle (2, 2);
        \draw[line width=0.5mm,  black] (2, 2) rectangle (3, 3);
        \draw[line width=0.5mm,  black] (3, 1) rectangle (4, 2);
        \draw[line width=0.5mm,  black] (2, 0) rectangle (3, 1);
        \draw[black] (0, 2) rectangle (2, 4);
        \draw[step=1cm, black] (0, 0) grid (2, 2);
    \end{tikzpicture}
    }\qquad 
    \caption{Partitioning of $\mathcal{IL}_{*}(X)$ at level $3$ in $2$D.}
    \label{fig:HODLRdD_interaction_split}
\end{center}
\end{figure}

It is well established that NCA can be applied under the standard (strong) admissibility condition, wherein admissible cluster pairs are either far-field or well-separated \cite{zhao2019fast,bebendorf2012constructing}. In contrast, the hybrid weak admissibility condition is defined by enforcing our weak admissibility criterion \cref{eq:weak_admis}, i.e., \( \text{adm}_{*}(X, Y) = \text{true} \), across levels $1$ through $\ell$ $(\ell < L)$. Under this weak admissibility condition, far-field and vertex-sharing clusters are admissible. To apply NCA, we appropriately partition the interaction list. More precisely, from level $1$ to $\ell$ of the tree, we partition the interaction list $\mathcal{IL}_{*} \bkt{X}$ into the \emph{far-field} interaction list denoted as $\mathcal{IL}_*^{far} \bkt{X}$ and the \emph{vertex-sharing} interaction list denoted as $\mathcal{IL}_*^{ver} \bkt{X}$, i.e., $\mathcal{IL}_{*} \bkt{X} = \mathcal{IL}_*^{far} \bkt{X} \bigcup \mathcal{IL}_*^{ver} \bkt{X}$ with $\mathcal{IL}_*^{far} \bkt{X} \bigcap \mathcal{IL}_*^{ver} \bkt{X} = \phi$. The visual representation of this partitioning is illustrated in \Cref{fig:HODLRdD_interaction_split}. 

Since the far-field clusters are separated, NCA with bottom to top pivot selection (B$2$T NCA) \cite{zhao2019fast, gujjula2022new} is applied to $\mathcal{IL}_*^{far}(X)$, and the corresponding bases are formed in a nested manner. On the other hand, ACA is applied to $\mathcal{IL}_*^{ver}(X)$, with bases constructed in a non-nested manner. Under the hybrid admissibility condition, adjacent clusters are treated as admissible at finer levels, and the corresponding matrix blocks are compressed using ACA.

The $\mathcal{H}_{hw}^{\prime}$-matrix for fast matrix-vector product is discussed in the following steps:
\begin{enumerate}
    \item Initialization of $\mathcal{H}_{hw}^{\prime}$-matrix representation (\Cref{subsec:H1_init})
    \item $\mathcal{H}_{hw}^{\prime}$-matrix-vector product (MVP) (\Cref{subsec:H1_mvp})
\end{enumerate}

\subsubsection{Initialization of \texorpdfstring{$\mathcal{H}_{hw}^{\prime}$}{H1}-matrix representation} \label{subsec:H1_init} We describe the initialization procedure in detail as follows:

\begin{enumerate}
    \item \underline{\emph{B$2$T NCA for far-field interaction.}} We apply the B$2$T NCA only to the far-field interactions, i.e., $\mathcal{IL}_{far}(X)$.  In the B$2$T NCA, the global row and column indices of the leaf clusters are considered to apply ACA and extract the pivots. For non-leaf clusters, ACA is applied to the pivots from the children-level clusters. The search spaces of the pivots for a particular cluster are illustrated in \Cref{fig:B2T_pivot_nHODLR}. This process recursively goes from bottom to top of the $2^d$ tree (from level $\ell$ to level $2$ of the tree) to generate all the required pivots. The B$2$T NCA is performed from level $\ell$ to level $2$ of the tree.  The pseudocode of B$2$T NCA is given in \Cref{algo:b2t_nca}.
    
    The translation operators corresponding to far-field interactions are constructed as described in \Cref{sec:nca}. We get the following sets of operators: 
 \begin{enumerate}
     \item P$2$M $\bkt{\mathcal{V}^{X}_{far}}$ and M$2$M $\bkt{\mathcal{\Tilde{V}}^{X X_c}_{far}}$, \quad $X_c \in child \bkt{X}$.
     \item M$2$L $\bkt{\mathcal{K}^{X,Y}_{far}}, \quad Y \in \mathcal{IL}_*^{far} \bkt{X}$.
     \item L$2$L $\bkt{\mathcal{\Tilde{U}}^{X_c X}_{far}}$, \quad $X \in parent \bkt{X_c}$ and L$2$P $\bkt{\mathcal{U}^{X}_{far}}$.
 \end{enumerate}

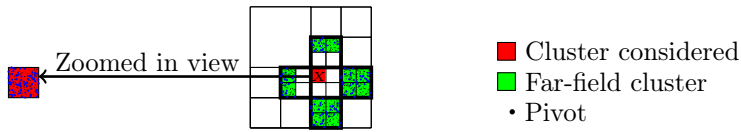
\begin{figure}[H]
\begin{center}
\captionsetup[subfloat]{labelformat=empty}
    \subfloat[]{
    \begin{tikzpicture}[scale=0.4]
        %%SELF
        \draw[draw=black,fill=red] (-8,1) rectangle (-7,2);
        \pgfmathsetseed{58678246}
        \foreach \x in {1,...,200}
        {
		\pgfmathsetmacro{\x}{rand/2 -7.5}
		\pgfmathsetmacro{\y}{rand/2 +1.5}
		\fill[black]    (\x,\y) circle (0.015);
        };
        \pgfmathsetseed{58678246}
        \foreach \x in {1,...,80}
        {
		\pgfmathsetmacro{\x}{rand/2 -7.5}
		\pgfmathsetmacro{\y}{rand/2 +1.5}
		\fill[blue]    (\x,\y) circle (0.035);
        };
        %%=================%%
        \draw[black] (0, 0) grid (4, 4);
        % \fill [green!50] (0,0) rectangle (1,1);
        % \fill [green!50] (0,1) rectangle (1,2);
        % \fill [green!50] (2,3) rectangle (3,4);
        % \fill [green!50] (3,3) rectangle (4,4);
        % \fill [cyan!50] (1,0) rectangle (2,1);
        % \fill [cyan] (1,2) rectangle (2,3);
        % \fill [cyan!50] (3,2) rectangle (4,3);
        % \fill [cyan!50] (3,0) rectangle (4,1);
        \fill [red] (2,1) rectangle (3,2);
        \fill [green] (1,1) rectangle (2,2);
        \fill [white] (1.5,1.5) rectangle (2,2);
        \fill [white] (2.5,1.5) rectangle (3,2);
        \fill [white] (2.5,1) rectangle (3,1.5);
        \fill [white] (2,1) rectangle (2.5,1.5);
        \fill [white] (2,2) rectangle (2.5,2.5);
        \fill [green] (2,2) rectangle (3,3);
        \fill [green] (2,0) rectangle (3,1);
        \fill [green] (3,1) rectangle (4,2);
        \node[anchor=north] at (2.27,2.15) {\tiny $X$};
        \node[anchor=north] at (1.5,2.8) {};
        \node[anchor=north] at (1.5,1.8) {};
        \node[anchor=north] at (.5,1.8) {};
        \draw[black] (0, 0) grid (4, 4);
        \draw[step=0.5cm, black] (3, 1) grid (4, 2);
        \draw[step=0.5cm, black] (2, 0) grid (3, 1);
        \draw[step=0.5cm, black] (2, 2) grid (3, 3);
        \draw[step=0.5cm, black] (2, 1) grid (3, 2);
        \draw[step=0.5cm, black] (1, 1) grid (2, 2);
        % \fill [cyan!25] (0,2) rectangle (2,4);
        \fill [white] (2,2) rectangle (2.5,2.5);
        \fill [white] (2.5,2) rectangle (3,2.5);
        \draw[line width=0.1mm,  black] (2.5, 2) rectangle (3, 2.5);
        \draw[line width=0.1mm,  black] (2, 2) rectangle (2.5, 2.5);
        \fill [white] (1.5,1) rectangle (2,1.5);
        \fill [white] (0,2) rectangle (2,4);
        \draw[line width=0.5mm,  black] (1, 1) rectangle (2, 2);
        \draw[line width=0.5mm,  black] (2, 2) rectangle (3, 3);
        \draw[line width=0.5mm,  black] (3, 1) rectangle (4, 2);
        \draw[line width=0.5mm,  black] (2, 0) rectangle (3, 1);
        \draw[black] (0, 2) rectangle (2, 4);
        \draw[step=1cm, black] (0, 0) grid (2, 2);
        \draw[step=0.5cm, black] (1, 1) grid (2, 2);

        \pgfmathsetseed{58678246}
        \foreach \x in {1,...,200}
        {
		\pgfmathsetmacro{\x}{rand/2+2.5}
		\pgfmathsetmacro{\y}{rand/2+0.5}
		\fill[black]    (\x,\y) circle (0.015);
        };

        \pgfmathsetseed{52213446}
        \foreach \x in {1,...,100}
        {
		\pgfmathsetmacro{\x}{rand/2+2.5}
		\pgfmathsetmacro{\y}{rand/4+2.75}
		\fill[black]    (\x,\y) circle (0.015);
        };

        \pgfmathsetseed{52213446}
        \foreach \x in {1,...,200}
        {
		\pgfmathsetmacro{\x}{rand/2+3.5}
		\pgfmathsetmacro{\y}{rand/2+1.5}
		\fill[black]    (\x,\y) circle (0.015);
        };

        \pgfmathsetseed{52213446}
        \foreach \x in {1,...,100}
        {
		\pgfmathsetmacro{\x}{rand/4+1.25}
		\pgfmathsetmacro{\y}{rand/2+1.5}
		\fill[black]    (\x,\y) circle (0.015);
        };

        \draw[->, line width=0.35mm] (2,1.7) -- (-7,1.7);
        \node at (-3.4,2.2) {Zoomed in view};
        %%%%%%%%%%%%%%%%%%%%%%%%%%%

        \pgfmathsetseed{58678246}
        \foreach \x in {1,...,40}
        {
		\pgfmathsetmacro{\x}{rand/2+2.5}
		\pgfmathsetmacro{\y}{rand/2+0.5}
		\fill[blue]    (\x,\y) circle (0.035);
        };

        \pgfmathsetseed{58678246}
        \foreach \x in {1,...,20}
        {
		\pgfmathsetmacro{\x}{rand/2+2.5}
		\pgfmathsetmacro{\y}{rand/4+2.75}
		\fill[blue]    (\x,\y) circle (0.035);
        };

        \pgfmathsetseed{58678246}
        \foreach \x in {1,...,40}
        {
		\pgfmathsetmacro{\x}{rand/2+3.5}
		\pgfmathsetmacro{\y}{rand/2+1.5}
		\fill[blue]    (\x,\y) circle (0.035);
        };

        \pgfmathsetseed{58678246}
        \foreach \x in {1,...,20}
        {
		\pgfmathsetmacro{\x}{rand/4+1.25}
		\pgfmathsetmacro{\y}{rand/2+1.5}
		\fill[blue]    (\x,\y) circle (0.035);
        };

    \end{tikzpicture}
    }\qquad \qquad
    \subfloat{
        \begin{tikzpicture}
            [%%%%%%%%%%%%%%%%%%%%%%%%%%%%%%
            box/.style={rectangle,draw=black, minimum size=0.2cm},scale=0.2
            ]%%%%%%%%%%%%%%%%%%%%%%%%%%%%%%
            \node[box,fill=red,,font=\tiny,label=right:Cluster considered,  anchor=west] at (5,8){};
            \node[box,fill=green,,font=\tiny,label=right:Far-field cluster,  anchor=west] at (5,6){};
            \node[circle,fill=black,inner sep=0pt,minimum size=2pt,label=right:{Pivot}] (a) at (6,4) {};
            % \node[box,fill=green,,font=\tiny,label=right:Far-field compressed block,  anchor=west] at (40,8){};
        \end{tikzpicture}
    }   
    \caption{Illustration of search spaces of the pivots in B$2$T NCA.}
    \label{fig:B2T_pivot_nHODLR}
\end{center}
\end{figure}

\begin{algorithm}[H]
 \small
\caption{B$2$T NCA for far-field interactions.}
\label{algo:b2t_nca}
\begin{algorithmic}[1]
    \For{\texttt{$l=\ell : -1 : 2$}} 
        \Comment{Traverse the tree from bottom to top direction}
				\For{\texttt{$i=1:2^{dl}$}}
                    \State $X \leftarrow $ $i^{th}$ cluster at level $l$ of tree.
                    \State $I.clear()$
                    \State $J.clear()$
                    \If{$l \neq \ell$} 
                        \State $I.insert \bkt{\bigcup \tau_{X_c}}$ \Comment{$X_c \in child(X)$, $c = 1 : 2^d$}
                        \For{$ Y \in $ \texttt{$\mathcal{IL}_*^{far}\bkt{{X}}$}}
						      \State $J.insert \bkt{\bigcup \tau_{Y_c}}$ \Comment{$Y_c \in child(Y)$, $c = 1 : 2^d$}
					    \EndFor
                    \Else
                    \State $I = X$ 
                    \For{ $Y \in $  \texttt{$\mathcal{IL}_*^{far}\bkt{{X}}$}}
						      \State $J.insert \bkt{Y}$ 
					\EndFor
                    \EndIf 
                \State $[\tau_X , \sigma_X] = \text{ACA} \bkt{I, J, \epsilon}$ \Comment{Perform ACA with tolerance $\epsilon$.}
				\EndFor
	\EndFor
\end{algorithmic}
\end{algorithm}

    \item  \underline{\textit{ACA for vertex-sharing interaction.}} The vertex-sharing interaction is compressed using partially pivoted ACA \cite{aca}. For each admissible block $X \times Y$, $Y \in \mathcal{IL}_*^{ver} \bkt{X}$, the vertex-sharing interaction block matrix $K_{X, Y}$ is compressed by partially pivoted ACA with user-given tolerance $\epsilon$ as follows:
    \begin{align*}
        K_{X, Y} \approx U V^*.
    \end{align*}
    
\Cref{algo:pure_aca_ver} outlines the vertex-sharing compression routine.

\begin{algorithm}[H]
 \small
	\caption{ACA for vertex-sharing interactions.}\label{algo:pure_aca_ver}
	\begin{algorithmic}[1]
			\For{\texttt{$l=1:\ell$}}  \Comment{Traverse the tree from level $1$ to $\ell$}
				\For{\texttt{$i=1:2^{dl}$}}
                    \State $X \leftarrow $ $i^{th}$ cluster at level $l$ of tree.
					% \State $I = t^X$ \Comment{$t^X = $ \text{ global index set of } $X$}
					\For{$Y \in \mathcal{IL}_*^{ver} \bkt{X}$}
						% \State $J = s^Y$ \Comment{$s^Y = $ \text{ global index set of } $Y$}
                          \State $[U, V] = \text{ACA}(X, Y, \epsilon)$ \Comment{Perform ACA to $K_{X,Y}$ with tolerance $\epsilon$.}
					\EndFor
				\EndFor
			\EndFor
	\end{algorithmic}
\end{algorithm}

\item \underline{\textit{ACA at the switching level.}} At the switching level $\ell$, the block matrix $K_{X,Y}$ is compressed using partially pivoted ACA, where $X \in \mathcal{T}^{\ell}$ and $Y \in \mathcal{N}_* \bkt{X}$. The detailed algorithm is presented in \Cref{algo:joining_level_H1}.

\begin{algorithm}[H]
 \small
\caption{ACA at the switching level $\ell$ $(\ell < L)$.}
\label{algo:joining_level_H1}
\begin{algorithmic}[1]
				\For{\texttt{$i=1:2^{d \ell}$}}    \Comment{At switching level $\ell$}
                    \State $X \leftarrow $ $i^{th}$ cluster at the level $\ell$ of tree.
					% \State $I =t^{X}$ \Comment{$t^X =$ \text{global index set of } $X$}
					\For{ $Y \in \mathcal{N}_* \bkt{X}$}
						% \State $J = s^Y$ \Comment{$s^Y =$ \text{global index set of } $Y$}
                        \State $[U , V] = \text{ACA} \bkt{X, Y, \epsilon}$ \Comment{Perform ACA to $K_{X,Y}$ with tolerance $\epsilon$.}
					\EndFor
				\EndFor
\end{algorithmic}
\end{algorithm}

\item \underline{\textit{ACA at finer levels.}} The finer level block matrices from level $\ell+1$ to level $L$, which are relatively smaller in size, are compressed using ACA as described in \Cref{algo:hodlr_finer_H1}.

\begin{algorithm}[H]
 \small
\caption{ACA at level $\ell +1$ to level $L$.}
\label{algo:hodlr_finer_H1}
\begin{algorithmic}[1]
    \For{\texttt{$l=\ell + 1:L$}} 
        \Comment{Traverse the tree from level $\ell +1$ to $L$}
				\For{\texttt{$i=1:2^{dl}$}}
                    \State $X \leftarrow $ $i^{th}$ cluster at level $l$ of tree.
					% \State $I =t^{X}$ \Comment{$t^X =$ \text{global index set of } $X$}
					\For{ $Y \in $  \texttt{$siblings\bkt{{X}}$}}
						% \State $J = s^Y$ \Comment{$s^Y =$ \text{global index set of } $Y$}
                        \State $[U , V] = \text{ACA} \bkt{X, Y, \epsilon}$ \Comment{Perform ACA to $K_{X,Y}$ with tolerance $\epsilon$.}
					\EndFor
				\EndFor
	\EndFor
\end{algorithmic}
\end{algorithm}
 
\end{enumerate}

\subsubsection{\texorpdfstring{$\mathcal{H}_{hw}^{\prime}$}{H1}-matrix-vector product (MVP)} \label{subsec:H1_mvp} 
Once we have the $\mathcal{H}_{hw}^{\prime}$-matrix representation, it can be applied to a column vector $\pmb{q}$ to compute fast MVP or potential, i.e., $\pmb{\Tilde{\phi}} = \Tilde{K} \pmb{q}$. Both nested and non-nested approaches are used to compute MVP. The detailed procedure is given below. 

\begin{enumerate}
\item \underline{\textit{Far-field MVP.}} We compute the MVP associated with far-field interaction using upward, transverse and downward tree traversal, similar to the $\mathcal{H}^2$-MVP algorithm.
\begin{enumerate}
    \item Upward traversal
        \begin{itemize}
            \item Particles to multipole (P$2$M) at level $\ell$: \\For all clusters $X \in \mathcal{T}^{\ell}$, calculate \\ $v^{far}_{X} = \mathcal{V}_{X}^{far^*} q_X$.
            \item Multipole to multipole (M$2$M) at levels $\ell-1 \geq l \geq 2$: \\For all clusters $X \in \mathcal{T}^l$, calculate \\ $v^{far}_{X} = \dsum_{X_c \in child \bkt{X}}  \Tilde{\mathcal{V}}_{X X_c}^{far^*} v_{X_c}^{far}$. 
        \end{itemize}
    \item  Transverse traversal
        \begin{itemize}
            \item Multipole to local (M$2$L) at levels $2 \leq l \leq \ell$: \\For all cluster $X \in \mathcal{T}^l$, calculate \\
             $u^{far}_{X} = \dsum_{Y \in \mathcal{IL}_*^{far} \bkt{X}}  \mathcal{K}_{X,Y}^{far} v_{Y}^{far}$.
        \end{itemize}
    \item  Downward traversal
        \begin{itemize}
            \item Local to local (L$2$L) at levels $2 \leq l \leq \ell-1$: \\For all clusters $X \in \mathcal{T}^l$, calculate \\
             $u^{far}_{X_c} \gets u^{far}_{X_c} +    \Tilde{\mathcal{U}}_{X_c X}^{far} u^{far}_{X} , \quad X \in parent \bkt{X_c} $.
            \item Local to particles (L$2$P) at level $\ell$: \\For all leaf clusters $X \in \mathcal{T}^{\ell}$, calculate \\ $\phi^{far}_{X} = \mathcal{U}_{X}^{far} u^{far}_{X}$.
        \end{itemize}
\end{enumerate}

 \item \underline{\textit{Vertex-sharing MVP.}}   The MVP associated with the vertex-sharing interaction is calculated through a single tree traversal in a non-nested manner, analogous to the $\mathcal{H}$-MVP. The details are presented in \Cref{algo:H1_ver}.

    % \textbf{Pseudocode of calculating the vertex-sharing potential.} The pseudocode outlining the calculation of the vertex-sharing potential ($\phi_{ver}$) is presented here.
 \begin{algorithm}[H]
  \small
	\caption{Vertex-sharing \mvp}\label{algo:H1_ver}
	\begin{algorithmic}[1]
	\State $\phi^{ver} \gets 0 \in \mathbb{R}^{N \times 1}$
			\For{\texttt{$l=1:\ell$}} 
				\For{\texttt{$i=1:2^{dl}$}}
                    \State $X \leftarrow $ $i^{th}$ cluster at level $l$ of tree.
					\For{$Y \in \mathcal{IL}_*^{ver} \bkt{X}$}. 
                        \State $\phi^{ver}_X \gets \phi^{ver}_X + U \bkt{V^* q_Y}$ \Comment{Low-rank MVP}
					\EndFor
				\EndFor
			\EndFor
	\end{algorithmic}
\end{algorithm}

\item \underline{\textit{Switching level MVP.}} At switching level $\ell$ $(\ell < L)$, the MVP is computed in a non-nested manner as described in \Cref{algo:joining_level}.

\begin{algorithm}[H]
 \small
\caption{Switching level \mvp}
\label{algo:joining_level}
\begin{algorithmic}[1]
\State $\phi_X^{\prime} \gets 0 \in \mathbb{R}^{N \times 1}$
				\For{\texttt{$i=1:2^{d \ell}$}}    
                    \State $X \leftarrow $ $i^{th}$ cluster at the level $\ell$ of tree.
					% \State $I =t^{X}$ \Comment{$t^X =$ \text{global index set of } $X$}
					\For{ $Y \in \mathcal{N}_* \bkt{X}$}
						% \State $J = s^Y$ \Comment{$s^Y =$ \text{global index set of } $Y$}
                        \State $\phi_X^{\prime} \gets \phi_X^{\prime} + U \bkt{V^* q_Y}$ \Comment{Low-rank MVP}
					\EndFor
				\EndFor
\end{algorithmic}
\end{algorithm}

\item \underline{\textit{Finer levels MVP.}} At the finer levels, from $\ell+1$ to $L$, the MVP is computed in a non-nested manner, following the procedure outlined in \Cref{algo:hodlr_finer}.

\begin{algorithm}[H]
 \small
\caption{Finer levels \mvp}
\label{algo:hodlr_finer}
\begin{algorithmic}[1]
    \For{\texttt{$l=\ell + 1:L$}} 
				\For{\texttt{$i=1:2^{dl}$}}
                    \State $X \leftarrow $ $i^{th}$ cluster at level $l$ of tree.
					% \State $I =t^{X}$ \Comment{$t^X =$ \text{global index set of } $X$}
					\For{ $Y \in $  \texttt{$siblings\bkt{{X}}$}}
						% \State $J = s^Y$ \Comment{$s^Y =$ \text{global index set of } $Y$}
                        \State $\phi_X^{\prime} \gets \phi_X^{\prime} + U \bkt{V^* q_Y}$ \Comment{Low-rank MVP}
					\EndFor
				\EndFor
	\EndFor
\end{algorithmic}
\end{algorithm}
\end{enumerate}

\subsubsection*{Final computed MVP}
            For every leaf cluster $X \in \mathcal{T}^{L}$, we add the MVP associated with self-interaction, which is a direct computation $\bkt{K_{X,X} q_{X}}$. Hence, the final computed MVP of the leaf cluster $X$ is given by
    \begin{equation*}
        \Tilde{\phi}_{X} = \underbrace{\phi_{X}^{far}}_{ \text{\tiny \makecell{Far-field MVP\\ (Nested)}}} + \underbrace{\phi_{X}^{ver}}_{\text{\tiny \makecell{Vertex-sharing MVP\\ (Non-nested)}}} + \underbrace{\phi^{\prime}_X}_{ \text{\tiny \makecell{Finer levels MVP\\ (Non-nested)}}} +  \underbrace{K_{X,X} q_{X}}_{\text{\tiny \makecell{Self interaction MVP\\ (Direct)}}}.
    \end{equation*}

Since $\Tilde{\phi}_{X_i}$ denotes the MVP corresponding to the leaf cluster $X_i$, the computed MVP can be expressed in MATLAB notation as $\pmb{\Tilde{\phi}} =\big [\Tilde{\phi}_{X_1};\Tilde{\phi}_{X_2}; \cdots ;\Tilde{\phi}_{X_{2^{d L}}} \big ]$.

\begin{remark}
    For $\mathcal{H}_{hw}^{\prime}$-matrix format, we consider $\ell<L$. If $\ell=L$, the neighboring blocks are not compressed and are instead treated as usual dense blocks. In this case, the resulting matrix structure reduces to the $\bkt{\mathcal{H}^2 + \mathcal{H}}_*$-matrix format of \cite{nHODLRdD}. On the other hand, $\ell=0$ yields the HODLR matrix format ($2^d$-tree-based), since in this case only \Cref{algo:hodlr_finer_H1} and \Cref{algo:hodlr_finer} remain active.  
\end{remark}

\subsubsection{Complexity analysis of \texorpdfstring{$\mathcal{H}_{hw}^{\prime}$-matrix}{H1}}
The cost of B$2$T NCA is $\mathcal{O}\bkt{N}$ \cite{zhao2019fast}, while the cost of ACA $\mathcal{O}\bkt{N \log(N)}$. Therefore, the overall initialization cost scales quasi-linearly.

The cost of far-field MVP is $\mathcal{O}\bkt{N}$ and all other steps scale quasi-linearly.

\subsection{\texorpdfstring{$\mathcal{H}_{hw}^{\prime \prime}$}{H1}-matrix} \label{subsec:H2}
The $\mathcal{H}_{hw}^{\prime}$-matrix algorithm can be further improved by employing nested bases for vertex-sharing interactions, resulting in the $\mathcal{H}_{hw}^{\prime\prime}$-matrix algorithm. In $\mathcal{H}_{hw}^{\prime\prime}$-matrix representation, vertex-sharing interactions are compressed using NCA with top to bottom pivot selection (T$2$B NCA) \cite{zhao2019fast}, leading to improved memory requirements and MVP computations. The difference between the $\mathcal{H}_{hw}^{\prime}$ and the $\mathcal{H}_{hw}^{\prime \prime}$ representations is the vertex-sharing interactions compression routine.

The $\mathcal{H}_{hw}^{\prime \prime}$-matrix for fast matrix-vector product is discussed in the following steps:
\begin{enumerate}
    \item Initialization of $\mathcal{H}_{hw}^{\prime \prime}$-matrix representation (\Cref{subsec:H2_init})
    \item $\mathcal{H}_{hw}^{\prime \prime}$-matrix-vector product (MVP) (\Cref{subsec:H2_mvp})
\end{enumerate}

\subsubsection{Initialization of \texorpdfstring{$\mathcal{H}_{hw}^{\prime \prime}$}{Hh1}-matrix representation} \label{subsec:H2_init} The initialization procedure differs from that of the $\mathcal{H}_{hw}^{\prime}$-matrix only in Step $2$ and is detailed as follows:

\begin{enumerate}
    \item \underline{\emph{B$2$T NCA for far-field interaction.}} Analogous to \Cref{algo:b2t_nca}.
    
    \item  \underline{\textit{T$2$B NCA for vertex-sharing interaction.}} The T$2$B NCA is applied only to the vertex-sharing interactions, i.e., $\mathcal{IL}_{ver} \bkt{X}$. The parent-level pivots are added to the child clusters, and the tree is traversed in a top to bottom manner (from level $1$ to level $\ell$ of the tree), as detailed in \Cref{algo:t2b_nca}. The search spaces of the pivots for a particular cluster are illustrated in \Cref{fig:T2B_pivot_nHODLR}.

\begin{figure}[H]
\begin{center}
\captionsetup[subfloat]{labelformat=empty}
    \subfloat[]{
    \begin{tikzpicture}[scale=0.4]
        %%SELF
        \draw[draw=black,fill=red] (-8,1) rectangle (-7,2);
        \pgfmathsetseed{58678246}
        \foreach \x in {1,...,200}
        {
		\pgfmathsetmacro{\x}{rand/2 -7.5}
		\pgfmathsetmacro{\y}{rand/2 +1.5}
		\fill[black]    (\x,\y) circle (0.015);
        };
        \pgfmathsetseed{58678246}
        \foreach \x in {1,...,80}
        {
		\pgfmathsetmacro{\x}{rand/2 -7.5}
		\pgfmathsetmacro{\y}{rand/2 +1.5}
		\fill[blue]    (\x,\y) circle (0.035);
        };
        %%=================%%
        \draw[black] (0, 0) grid (4, 4);
        % \fill [green!50] (0,0) rectangle (1,1);
        % \fill [green!50] (0,1) rectangle (1,2);
        % \fill [green!50] (2,3) rectangle (3,4);
        % \fill [green!50] (3,3) rectangle (4,4);
        \fill [cyan!25] (1,0) rectangle (2,1);
        % \fill [cyan] (1,2) rectangle (2,3);
        \fill [cyan!25] (3,2) rectangle (4,3);
        \fill [cyan!25] (3,0) rectangle (4,1);
        \fill [red] (2,1.5) rectangle (2.5,2);
        % \fill [green] (1,1) rectangle (2,2);
        % \fill [orange] (1.5,1.5) rectangle (2,2);
        % \fill [orange] (2.5,1.5) rectangle (3,2);
        \fill [cyan] (2.5,1) rectangle (3,1.5);
        % \fill [orange] (2,1) rectangle (2.5,1.5);
        % \fill [orange] (2,2) rectangle (2.5,2.5);
        % \fill [green] (2,2) rectangle (3,3);
        % \fill [green] (2,0) rectangle (3,1);
        % \fill [green] (3,1) rectangle (4,2);
        \node[anchor=north] at (1.5,2.8) {};
        \node[anchor=north] at (1.5,1.8) {};
        \node[anchor=north] at (.5,1.8) {};
        \draw[black] (0, 0) grid (4, 4);
        \draw[step=0.5cm, black] (3, 1) grid (4, 2);
        \draw[step=0.5cm, black] (2, 0) grid (3, 1);
        \draw[step=0.5cm, black] (2, 2) grid (3, 3);
        \draw[step=0.5cm, black] (2, 1) grid (3, 2);
        \draw[step=0.5cm, black] (1, 1) grid (2, 2);
        \fill [white] (0,2) rectangle (2,4);
        % \fill [orange] (2,2) rectangle (2.5,2.5);
        \fill [cyan] (2.5,2) rectangle (3,2.5);
        \fill [cyan] (1.5,1) rectangle (2,1.5);
        \draw[line width=0.5mm,  black] (1, 1) rectangle (2, 2);
        \draw[line width=0.5mm,  black] (2, 2) rectangle (3, 3);
        \draw[line width=0.5mm,  black] (3, 1) rectangle (4, 2);
        \draw[line width=0.5mm,  black] (2, 0) rectangle (3, 1);
        \draw[black] (0, 2) rectangle (2, 4);
        \draw[step=1cm, black] (0, 0) grid (2, 2);
        \node[anchor=north] at (2.27,2.15) {\tiny $X$};

        \pgfmathsetseed{58678246}
        \foreach \x in {1,...,200}
        {
		\pgfmathsetmacro{\x}{rand/2+3.5}
		\pgfmathsetmacro{\y}{rand/2+0.5}
		\fill[black]    (\x,\y) circle (0.015);
        };

        \pgfmathsetseed{52213446}
        \foreach \x in {1,...,200}
        {
		\pgfmathsetmacro{\x}{rand/2+3.5}
		\pgfmathsetmacro{\y}{rand/2+2.5}
		\fill[black]    (\x,\y) circle (0.015);
        };

        \pgfmathsetseed{52213446}
        \foreach \x in {1,...,200}
        {
		\pgfmathsetmacro{\x}{rand/2+1.5}
		\pgfmathsetmacro{\y}{rand/2+0.5}
		\fill[black]    (\x,\y) circle (0.015);
        };

        \pgfmathsetseed{58678248}
        \foreach \x in {1,...,80}
        {
		\pgfmathsetmacro{\x}{rand/4+1.75}
		\pgfmathsetmacro{\y}{rand/4+1.25}
		\fill[blue]    (\x,\y) circle (0.015);
        };

        \pgfmathsetseed{58678249}
        \foreach \x in {1,...,80}
        {
		\pgfmathsetmacro{\x}{rand/4+2.75}
		\pgfmathsetmacro{\y}{rand/4+1.25}
		\fill[blue]    (\x,\y) circle (0.015);
        };

        \pgfmathsetseed{58678257}
        \foreach \x in {1,...,80}
        {
		\pgfmathsetmacro{\x}{rand/4+2.75}
		\pgfmathsetmacro{\y}{rand/4+2.25}
		\fill[blue]    (\x,\y) circle (0.015);
        };

        \draw[->, line width=0.35mm] (2,1.7) -- (-7,1.7);
        \node at (-3.4,2.2) {Zoomed in view};
        %%%%%%%%%%%%%%%%%%%%%%%%%%%

        \pgfmathsetseed{58678246}
        \foreach \x in {1,...,40}
        {
		\pgfmathsetmacro{\x}{rand/2+3.5}
		\pgfmathsetmacro{\y}{rand/2+0.5}
		\fill[blue]    (\x,\y) circle (0.035);
        };

        \pgfmathsetseed{58678246}
        \foreach \x in {1,...,40}
        {
		\pgfmathsetmacro{\x}{rand/2+3.5}
		\pgfmathsetmacro{\y}{rand/2+2.5}
		\fill[blue]    (\x,\y) circle (0.035);
        };

        \pgfmathsetseed{58678246}
        \foreach \x in {1,...,40}
        {
		\pgfmathsetmacro{\x}{rand/2+1.5}
		\pgfmathsetmacro{\y}{rand/2+0.5}
		\fill[blue]    (\x,\y) circle (0.035);
        };

        \pgfmathsetseed{58678248}
        \foreach \x in {1,...,20}
        {
		\pgfmathsetmacro{\x}{rand/4+1.75}
		\pgfmathsetmacro{\y}{rand/4+1.25}
		\fill[blue]    (\x,\y) circle (0.035);
        };

        \pgfmathsetseed{58678249}
        \foreach \x in {1,...,20}
        {
		\pgfmathsetmacro{\x}{rand/4+2.75}
		\pgfmathsetmacro{\y}{rand/4+1.25}
		\fill[blue]    (\x,\y) circle (0.035);
        };

        \pgfmathsetseed{58678257}
        \foreach \x in {1,...,20}
        {
		\pgfmathsetmacro{\x}{rand/4+2.75}
		\pgfmathsetmacro{\y}{rand/4+2.25}
		\fill[blue]    (\x,\y) circle (0.035);
        };

    \end{tikzpicture}
    }\qquad \qquad
    \subfloat{
        \begin{tikzpicture}
            [%%%%%%%%%%%%%%%%%%%%%%%%%%%%%%
            box/.style={rectangle,draw=black, minimum size=0.2cm},scale=0.2
            ]%%%%%%%%%%%%%%%%%%%%%%%%%%%%%%
            \node[box,fill=red,,font=\tiny,label=right:Cluster considered,  anchor=west] at (5,8){};
            \node[box,fill=cyan!25,,font=\tiny,label=right:Parent-level vertex-sharing cluster ,  anchor=west] at (5,6){};
            \node[box,fill=cyan,,font=\tiny,label=right:Vertex-sharing cluster,  anchor=west] at (5,4){};
            \node[circle,fill=black,inner sep=0pt,minimum size=2pt,label=right:{Pivot}] (a) at (6,2) {};
            % \node[box,fill=green,,font=\tiny,label=right:Far-field compressed block,  anchor=west] at (40,8){};
        \end{tikzpicture}
    }   
    \caption{Illustration of search spaces of the pivots in T$2$B NCA for $\mathcal{H}_{hw}^{\prime \prime}$-matrix.}
    \label{fig:T2B_pivot_nHODLR}
\end{center}
\end{figure}
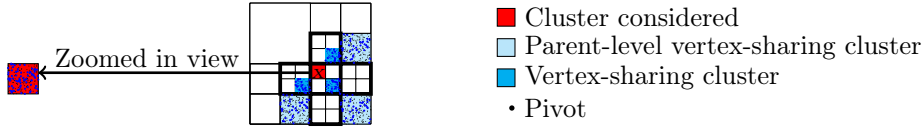

\begin{algorithm}[H]
 \small
\caption{T$2$B NCA for vertex-sharing interactions.}
\label{algo:t2b_nca}
\begin{algorithmic}[1]
    \For{\texttt{$l=1:\ell$}} 
        \Comment{Traverse the tree from top to bottom direction}
				\For{\texttt{$i=1:2^{dl}$}}
                    \State $X \leftarrow $ $i^{th}$ cluster at level $l$ of tree.
					\State $I=X$
					\For{ $Y \in $  \texttt{$\mathcal{IL}_{ver}\bkt{{X}}$}}
						\State $J=Y$
                        \If{$l > 1$} 
                            \State $J.insert \bkt{\sigma_{parent\bkt{X}}}$  \Comment{Add parent-level pivots}
                        \EndIf 
					\EndFor
                \State $[\tau_X , \sigma_X] = \text{ACA} \bkt{I, J, \epsilon}$ \Comment{Perform ACA with tolerance $\epsilon$}
				\EndFor
	\EndFor
\end{algorithmic}
\end{algorithm}

    The translation operators corresponding to vertex-sharing interactions are constructed as described in \Cref{sec:nca}. We get the following sets of operators: 
 \begin{enumerate}
     \item P$2$M $\bkt{\mathcal{V}_{X}^{ver}}$ and M$2$M $\bkt{\mathcal{\Tilde{V}}_{X X_c}^{ver}}$, \quad $X_c \in child \bkt{X}$.
     \item M$2$L $\bkt{\mathcal{K}_{X,Y}^{ver}}, \quad Y \in \mathcal{IL}_*^{ver} \bkt{X}$.
     \item L$2$L $\bkt{\mathcal{\Tilde{U}}_{X_c X}^{ver}}$, \quad $X \in parent \bkt{X_c}$ and L$2$P $\bkt{\mathcal{U}_{X}^{ver}}$.
 \end{enumerate}

\item \underline{\textit{ACA at the switching level.}} Analogous to \Cref{algo:joining_level_H1}.

\item \underline{\textit{ACA at finer levels.}} Analogous to \Cref{algo:hodlr_finer_H1}.
 
\end{enumerate}

\subsubsection{\texorpdfstring{$\mathcal{H}_{hw}^{\prime \prime}$}{H1}-matrix-vector product (MVP)} \label{subsec:H2_mvp} 
We apply $\mathcal{H}_{hw}^{\prime \prime}$-matrix representation to a column vector $\pmb{q}$ to compute fast MVP, i.e., $\pmb{\Tilde{\phi}} = \Tilde{K} \pmb{q}$. The MVPs associated with the far-field and vertex-sharing are evaluated in a nested manner, while MVPs at the finer levels are obtained in a non-nested manner, as follows:

\begin{enumerate}
\item \underline{\textit{Far-field MVP.}} Analogous to the corresponding step of $\mathcal{H}_{hw}^{\prime}$-matrix.

\item \underline{\textit{Vertex-sharing MVP}} The MVP corresponding to the vertex-sharing interactions is computed using upward, transverse and downward tree traversal, as given below.
\begin{enumerate}
    \item Upward traversal
        \begin{itemize}
            \item Particles to multipole (P$2$M) at level $\ell$: \\For all clusters $X \in \mathcal{T}^{\ell}$, calculate \\ $v^{ver}_{X} = \mathcal{V}_{X}^{ver^*} q_X$.
            \item Multipole to multipole (M$2$M) at levels $\ell-1 \geq l \geq 1$: \\For all clusters $X \in \mathcal{T}^l$, calculate \\ $v^{ver}_{X} = \dsum_{X_c \in child \bkt{X}}  \Tilde{\mathcal{V}}_{X X_c}^{ver^*} v_{X_c}^{ver}$. 
        \end{itemize}
    \item  Transverse traversal
        \begin{itemize}
            \item Multipole to local (M$2$L) at levels $1 \leq l \leq \ell$: \\For all cluster $X \in \mathcal{T}^l$, calculate \\
             $u^{ver}_{X} = \dsum_{Y \in \mathcal{IL}_*^{ver} \bkt{X}}  \mathcal{K}_{X,Y}^{ver} v_{Y}^{ver}$.
        \end{itemize}
    \item  Downward traversal
        \begin{itemize}
            \item Local to local (L$2$L) at levels $1 \leq l \leq \ell-1$: \\For all clusters $X \in \mathcal{T}^l$, calculate \\
             $u^{ver}_{X_c} \gets u^{ver}_{X_c} +    \Tilde{\mathcal{U}}_{X_c X}^{ver} u^{ver}_{X} , \quad X \in parent \bkt{X_c} $.
            \item Local to particles (L$2$P) at level $\ell$: \\For all leaf clusters $X \in \mathcal{T}^{\ell}$, calculate \\ $\phi^{ver}_{X} = \mathcal{U}_{X}^{ver} u^{ver}_{X}$.
        \end{itemize}
\end{enumerate}

\item \underline{\textit{Switching level MVP.}} Analogous to \Cref{algo:joining_level}.

\item \underline{\textit{Finer levels MVP.}} Analogous to \Cref{algo:hodlr_finer}.
\end{enumerate}

\subsubsection*{Final computed MVP}
            For every leaf cluster $X \in \mathcal{T}^{L}$, we add the MVP associated with self-interaction, which is a direct computation. Hence, the final computed MVP of the leaf cluster $X$ is given by
    \begin{equation*}
        \Tilde{\phi}_{X} = \underbrace{\phi_{X}^{far}}_{ \text{\tiny \makecell{Far-field MVP\\ (Nested)}}} + \underbrace{\phi_{X}^{ver}}_{\text{\tiny \makecell{Vertex-sharing MVP\\ (Nested)}}} + \underbrace{\phi^{\prime}_X}_{ \text{\tiny \makecell{Finer levels MVP\\ (Non-nested)}}} +  \underbrace{K_{X,X} q_{X}}_{\text{\tiny \makecell{Self interaction MVP\\ (Direct)}}}.
    \end{equation*}

The computed MVP can be expressed as $\pmb{\Tilde{\phi}} =\big [\Tilde{\phi}_{X_1};\Tilde{\phi}_{X_2}; \cdots ;\Tilde{\phi}_{X_{2^{d L}}} \big ]$.

\begin{remark}
    The B$2$T NCA is not suitable for vertex-sharing interactions because the rank associated with these interactions increases, as demonstrated in \cite{nHODLRdD}. Although the T$2$B NCA can be applied to the entire interaction list $\bkt{\mathcal{IL}_* (X)}$, doing so incurs a high initialization cost (see \cite{nHODLRdD}). Thus, a more efficient strategy is to partition the interaction list, applying the B$2$T NCA to far-field interactions and the T$2$B NCA to vertex-sharing interactions.
\end{remark}

\begin{remark}
    For $\mathcal{H}_{hw}^{\prime \prime}$-matrix format, we consider $\ell<L$. If $\ell=L$, the neighboring blocks are not compressed and are instead treated as usual dense blocks. In this case, the resulting matrix structure reduces to the $\mathcal{H}^2_*$-matrix format of \cite{nHODLRdD}. On the other hand, $\ell=0$ yields the HODLR matrix format ($2^d$-tree-based).  
\end{remark}

\subsubsection{Complexity analysis of \texorpdfstring{$\mathcal{H}_{hw}^{\prime \prime}$-matrix}{H2}}
The cost of T$2$B NCA (step $2$) is $\mathcal{O}\bkt{N \log(N)}$ \cite{zhao2019fast}, resulting in an overall quasi-linear complexity. Although the initialization of $\mathcal{H}_{hw}^{\prime \prime}$ and $\mathcal{H}_{hw}^{\prime}$ exhibits the same asymptotic complexity, $\mathcal{H}_{hw}^{\prime \prime}$ is observed to be slightly more computationally expensive because T$2$B NCA incurs a higher computational cost than ACA. However, the $\mathcal{H}_{hw}^{\prime \prime}$-MVP is faster than $\mathcal{H}_{hw}^{\prime}$-MVP due to the use of nested bases for vertex-sharing interactions. This is shown in \Cref{sec:num_results}.  

\subsection{Advantages of the proposed algorithms}
We investigate the condition under which the proposed hybrid weak admissibility-based algorithms offer advantages over the uniform weak admissibility-based (weak admissibility condition \cref{eq:weak_admis} is applied at all levels) algorithms proposed in \cite{nHODLRdD}.

Let $\mathcal{B}$, $\mathcal{N}$, and $\mathcal{D}$ denote the admissible, neighbor, and dense blocks, respectively. The subscript represents the admissibility condition, while the superscript indicates the number of blocks across levels of the tree. The front symbols $LR$ and $D$ represent the low-rank (cost: $\mathcal{O}\bkt{p(m+n)}$) and dense matrix (cost: $\mathcal{O}\bkt{mn}$) formats, respectively. More precisely, $LR\text{-}\mathcal{B}^{1:\ell}_*$ represents the storage cost of all admissible blocks stored in low-rank format under the weak admissibility condition \cref{eq:weak_admis} across the hierarchical levels $1$ through $\ell$.

The total storage cost of the proposed hybrid weak admissibility-based algorithms is given by
\begin{equation} \label{eq:s_hw}
    \mathcal{S}_{hw} = LR\text{-}\mathcal{B}^{1:\ell}_* + LR\text{-}\mathcal{N}^{\ell}_* + LR\text{-}\mathcal{B}^{\ell+1:L}_w + D\text{-}\mathcal{D}^L.
\end{equation}

The total storage cost of the (uniform) weak admissibility-based algorithms \cite{nHODLRdD} is given by
\begin{equation} \label{eq:s_*}
    \mathcal{S}_* = LR\text{-}\mathcal{B}^{1:L}_* + D\text{-}\mathcal{N}^{L}_* + D\text{-}\mathcal{D}^L = LR\text{-}\mathcal{B}^{1:\ell}_* + LR\text{-}\mathcal{B}^{\ell+1:L}_* + D\text{-}\mathcal{N}^{L}_* + D\text{-}\mathcal{D}^L.
\end{equation}

The proposed hybrid algorithms avoid constructing translation operators at the \emph{finer levels} of the tree, where the block sizes are relatively small. Instead, the corresponding matrix blocks are compressed in a non-nested manner, thereby reducing the computational overhead of constructing and applying translation operators. Furthermore, the proposed algorithms do not involve neighbor blocks, and the diagonal blocks are only dense blocks. This significantly reduces the number of dense matrix operations, resulting in improved computational efficiency and faster overall execution. Hence, the efficiency of the hybrid hierarchical matrix algorithms is determined by the cost incurred at the finer levels of the hierarchy. 

We establish the following condition that guarantees the computational advantage of the proposed hybrid algorithm.
From \cref{eq:s_hw} and \cref{eq:s_*}, we obtain
\begin{equation} \label{eq:comparison_eqn}
    \mathcal{S} \bkt{\ell} = \mathcal{S}_* - \mathcal{S}_{hw} = \bkt{LR\text{-}\mathcal{B}^{\ell+1:L}_* + D\text{-}\mathcal{N}^{L}_*} - \bkt{LR\text{-}\mathcal{N}^{\ell}_* + LR\text{-}\mathcal{B}^{\ell+1:L}_w}.
\end{equation}

Therefore, the proposed algorithms are advantageous in terms of storage requirements and MVP execution costs if $\mathcal{S} \bkt{\ell} >0$. We compare proposed hybrid algorithms with their counterparts based on (uniform) weak admissibility in \Cref{subsec:comparison}.

A natural question that arises is how to identify an optimal switching level $\ell^*$, which maximizes the performance gains over the uniformly weak admissibility-based counterparts, i.e.,
\begin{equation*} 
\ell^* = \arg \max_{0 < \ell < L} \bkt{\mathcal{S} \bkt{\ell}}. 
\end{equation*}

 In practice, one may evaluate $\mathcal{S}(\ell)$ for $\ell = L-1, L-2, \dots, 1$ to identify the optimal switching level. However, this can be computationally expensive. Therefore, in all numerical experiments reported in \Cref{sec:num_results}, we heuristically fix $\ell = L-1$ for simplicity.

\section{Numerical results} \label{sec:num_results}
%FLOPs, Improvement factor
In this section, we present numerical experiments to assess the performance of the proposed hybrid weak admissibility-based hierarchical matrices compared with their uniformly weak admissible counterparts. Furthermore, we demonstrate the applicability of the proposed algorithms to real-world problems, including fast principal component analysis in Gaussian processes and fast kernel ridge regression. In \Cref{tab:app_notations}, we summarize some notations used in this section.
  \begin{table}[H]
     \centering
     \resizebox{\textwidth}{!}{\begin{tabular}{|c|c|}
     \hline
    $n_{max}$ & \makecell{Maximum number of particles at leaf-clusters. We set $n_{max} = 125$.} \\
    \hline
    $\Tilde{K}$ & \makecell{$\Tilde{K}$ represent the hierarchical matrix representation of the original kernel matrix $K$.} \\
    \hline
    Mem. & \makecell{Total memory (in GB) needed to store the hierarchical low-rank representation ($\Tilde{K}$).} \\
    \hline
        $t_{init}$ &  \makecell{Time taken (in s) to initialize the hierarchical representation $\Tilde{K}$}  \\
    \hline
    $t_{FMVP}$  & \makecell{Time taken (in s) to compute the fast MVP $(\Tilde{K} \pmb{q})$.}   \\
    \hline
    $t_{DMVP}$  & \makecell{The time taken (in s) to compute the direct MVP.}   \\
    \hline
    ${RE}_{MVP}$ & \makecell{The $2$-norm relative error in the MVP.} \\
    \hline
    % Relative error in solution $({RE}_{sol})$ & \makecell{Let $\pmb{q}$ be a random column vector and the exact \mvp $\pmb{b} = K \pmb{q}$ (exact up to round-off). \\ We set $\pmb{b}$ as RHS and solve the system $K \pmb{\lambda} = \pmb{b}$ using the iterative solver, GMRES. \\ The $2$-norm relative error in the solution of the system is given by $\magn{\pmb{\lambda} - \pmb{q}}_2 / \magn{\pmb{q}}_2$ } \\
    % \hline
    % $\epsilon$ & \makecell{The NCA/ACA tolerance.} \\
    % \hline
    % GMRES stopping condition & \makecell{We use GMRES \cite{saad1986gmres} without restart, and GMRES routine will terminate if the relative residual is less than $\epsilon_{GMRES}$. \\ The \mvp part in the GMRES routine is accelerated using various \textbf{fast} MVP algorithms (\Cref{tab:nex_notation}). \\ We set $\epsilon_{GMRES} = 10^{-12}$ and $10^{-10}$ for the experiments in $2$D and $3$D, respectively.} \\
    \#iter & \makecell{The number of iterations of the CG to reach the stopping condition.} \\
    \hline
     \end{tabular}}
     \caption{Notations used in this section for various numerical experiments.}
     \label{tab:app_notations}
 \end{table}

 All the algorithms are implemented in \texttt{C++}, and the experiments are performed on an Intel Xeon(R) $2.30$GHz processor with $4$ \texttt{OpenMP} threads.

\subsection{Fast MVP in \texorpdfstring{$3$}{3}D} \label{subsec:comparison}
We consider $N$ uniformly distributed particles located at $\{\pmb{x}_i \}_{i=1}^N$ in the cube $[-1,1]^3$. The kernel function is chosen as the single-layer Laplacian in $3$D, and the $(i,j)^{th}$ entry of the kernel matrix $K \in \mathbb{R}^{N \times N}$ is given by
\begin{equation}
K(i,j) = \begin{cases}
            \dfrac{1}{r} = \dfrac{1}{\magn{\pmb{x}_i - \pmb{x}_j}_2} & \text{ if } i \neq j\\
             0 & \text{ otherwise}
            \end{cases}    
\end{equation}

We compare the proposed hybrid weak admissibility-based $\mathcal{H}_{hw}^{\prime}$ and $\mathcal{H}_{hw}^{\prime \prime}$ formats with their (uniform) weak admissibility-based counterparts, $\bkt{\mathcal{H}^2 + \mathcal{H}}_{*}$ and $\mathcal{H}^2_{*}$, respectively, where the weak admissibility condition \cref{eq:weak_admis} is applied to all the levels. Additionally, we compare these formats with the standard admissibility-based $\mathcal{H}^2_{\eta}$-matrix format (the interaction list is analogous to that in FMM when $\eta = \sqrt{d}$). It is worth noting that the leaf-level block size and the depth of the tree are kept the same for all $\mathcal{H}$-matrix formats. A summary of the different $\mathcal{H}$-matrix formats is provided in \Cref{tab:H_mat_summ} for ease of comparison.

 \begin{table}
     \resizebox{\textwidth}{!}{
    \begin{tabular}{|l|l|l|l|l|}
    \hline
     \textbf{$\mathcal{H}$-matrix} & \qquad \qquad \textbf{Short description} & \qquad \makecell{\textbf{Admissibility condition}/ \\ \textbf{Admissible clusters}} & \textbf{Compression technique(s)} &  \textbf{Type of bases} \\ \hline
     $\mathcal{H}_{hw}^{\prime}$        & This work, see \Cref{subsec:H1}.   & \makecell{Hybrid weak admissibility \cref{eq:hybrid_admiss}/ \\ \textbf{Coarser levels:} \emph{far-field}, \emph{vertex-sharing} \\ \textbf{Finer levels:} \emph{adjacent}} & \makecell{B$2$T NCA and ACA} & \makecell{Semi-nested and \\ non-nested}  \\ \hline
     $\mathcal{H}_{hw}^{\prime \prime}$        & This work, see \Cref{subsec:H2}.   & \makecell{Hybrid weak admissibility \cref{eq:hybrid_admiss}/ \\ \textbf{Coarser levels:} \emph{far-field}, \emph{vertex-sharing} \\ \textbf{Finer levels:} \emph{adjacent}} & \makecell{B$2$T NCA, T$2$B NCA, \\ and ACA} & \makecell{Nested and \\ non-nested}  \\ \hline
     $\bkt{\mathcal{H}^2 + \mathcal{H}}_{*}$             & \makecell{See $\bkt{\mathcal{H}^2 + \mathcal{H}}_{*}$-matrix of \cite{nHODLRdD}.}& \makecell{Weak admissibility condition \cref{eq:weak_admis}/ \\ \emph{far-field} and \emph{vertex-sharing}} & \quad B$2$T NCA and ACA & Semi-nested \\ \hline
     $\mathcal{H}^2_{*}$             & \makecell{See $\mathcal{H}^2_{*}$-matrix of \cite{nHODLRdD}.}& \makecell{Weak admissibility condition \cref{eq:weak_admis}/ \\ \emph{far-field} and \emph{vertex-sharing}} & \makecell{B$2$T NCA (far-field) and \\ T$2$B NCA (vertex-sharing)} & Nested \\ \hline
     $\mathcal{H}^2_{\eta}$             & \makecell{Standard admissibility-based $\mathcal{H}^2$-matrix.\\ We set $\eta = \sqrt{d}$ in \cref{eq:H_admiss}.}& \makecell{Standard admissibility condition \cref{eq:H_admiss}/ \\ \emph{far-field}} & \quad B$2$T NCA & Nested \\ \hline
    \end{tabular}}
     \caption{Different $\mathcal{H}$-matrix formats and their descriptions.}
     \label{tab:H_mat_summ}
 \end{table}

We randomly select $5$ different column vectors $(\pmb{q})$ and perform the $\mathcal{H}$-matrix-vector products $(\Tilde{K} \pmb{q})$ using different $\mathcal{H}$-matrix formats. The memory to store the $\mathcal{H}$-matrix format (Mem), average initialization time $(t_{init})$, $\mathcal{H}$-MVP time $(t_{FMVP})$, the mega FLOP count for $\mathcal{H}$-MVP (MFLOP), and the relative error in MVP $({RE}_{MVP})$ for different tolerances $(\epsilon)$ are tabulated in \Cref{table:smvp3d_lap_1,table:smvp3d_lap_2,table:smvp3d_lap_3}. We also report the direct MVP time $(t_{DMVP})$. 

\begin{table}
      \centering
      \resizebox{\textwidth}{!}{%
      \setlength\extrarowheight{0.9pt}
\begin{tabular}{|l|l|l|l|l|l|l|l|}
\hline
$N$                        & \makecell{$\mathcal{H}$-matrix} & Mem (GB) & $t_{init}$ (s) & $t_{FMVP}$ (s) & MFLOP & $t_{DMVP}$ (s)  & $RE_{MVP}$ \\ \hline
\multirow{5}{*}{64000}    & $\mathcal{H}_{hw}^{\prime}$ & 0.67               & 2.04   & 0.04           & 156.22    & 8.24 & 2.39E-03  \\ \cline{2-8} 
                         & $\mathcal{H}_{hw}^{\prime \prime}$  & \textbf{0.63}                & 2.36   & \textbf{0.03}           & \textbf{150.55}    & 8.24 & 2.36E-03  \\ \cline{2-8} 
                         & $\bkt{\mathcal{H}^2 + \mathcal{H}}_{*}$  & 0.75                & 1.77   & 0.06           & 175.02    & 8.24 & 1.53E-03  \\ \cline{2-8} 
                         & $\mathcal{H}^2_{*}$  & 0.72                & 2.55   & 0.05           & 171.43    & 8.24 & 1.59E-03  \\ \cline{2-8} 
                         & $\mathcal{H}^2_{\eta}$  & 0.84                & 2.19   & 0.08           & 202.73    & 8.24 & 1.42E-03  \\ \hline \hline
\multirow{5}{*}{512000}    & $\mathcal{H}_{hw}^{\prime}$ & 6.07               & 27.26   & 0.33           & 1415.69    & 527.61 & 2.24E-03  \\ \cline{2-8} 
                         & $\mathcal{H}_{hw}^{\prime \prime}$  & \textbf{5.36}                & 35.46   & \textbf{0.28}           & \textbf{1277.87}    & 527.61 & 2.27E-03  \\ \cline{2-8} 
                         & $\bkt{\mathcal{H}^2 + \mathcal{H}}_{*}$  & 7.41                & 21.15   & 0.52           & 1749.71    & 527.61 & 2.17E-03  \\ \cline{2-8} 
                         & $\mathcal{H}^2_{*}$  & 6.73                & 34.93   & 0.48           & 1618.54    & 527.61 & 2.28E-03  \\ \cline{2-8} 
                         & $\mathcal{H}^2_{\eta}$  & 8.32                & 32.53   & 0.72           & 2015.33    & 527.61 & 2.56E-03  \\ \hline \hline
\multirow{5}{*}{1000000}    & $\mathcal{H}_{hw}^{\prime}$ &    8.82      & 46.29    & 0.62          & 2072.14    & 2038.87 & 5.68E-03 \\ \cline{2-8} 
                         & $\mathcal{H}_{hw}^{\prime \prime}$  & \textbf{6.59}                & 128.24   & \textbf{0.57}           & \textbf{1617.63}    & 2038.87 & 5.67E-03 \\ \cline{2-8} 
                         & $\bkt{\mathcal{H}^2 + \mathcal{H}}_{*}$  & 14.66                & 64.18   & 1.12           & 3535.83    & 2038.87 & 1.19E-03  \\ \cline{2-8} 
                         & $\mathcal{H}^2_{*}$  & 12.82                & 154.17   & 1.05           & 3173.41    & 2038.87 & 1.22E-03  \\ \cline{2-8} 
                         & $\mathcal{H}^2_{\eta}$  & 15.46                & 90.06   & 1.27           & 3834.68    & 2038.87 & 1.21E-03  \\ \hline \hline
\multirow{5}{*}{2197000}    & $\mathcal{H}_{hw}^{\prime}$ & 21.10               & 130.22   & 1.24           & 4975.42    & 10062.11 & 1.72E-03 \\ \cline{2-8} 
                         & $\mathcal{H}_{hw}^{\prime \prime}$  & \textbf{16.84}               & 196.17   & \textbf{1.08}           & \textbf{4061.68}    & 10062.11 & 2.31E-03 \\ \cline{2-8} 
                         & $\bkt{\mathcal{H}^2 + \mathcal{H}}_{*}$  & 28.84                & 121.06   & 2.12           & 6911.54    & 10062.11 & 1.08E-03  \\ \cline{2-8} 
                         & $\mathcal{H}^2_{*}$  & 25.75                & 317.42   & 2.12           & 6288.31    & 10062.11 & 1.78E-03  \\ \cline{2-8} 
                         & $\mathcal{H}^2_{\eta}$  & 30.77                & 166.23   & 2.52           & 7544.71    & 10062.11 & 1.55E-03  \\ \hline 
\end{tabular}
}
    \caption{Comparison between different $\mathcal{H}$-matrix formats when $\epsilon = 10^{-3}$.}
    \label{table:smvp3d_lap_1}
\end{table}

\begin{table}
      \centering
      \resizebox{\textwidth}{!}{%
      \setlength\extrarowheight{0.9pt}
\begin{tabular}{|l|l|l|l|l|l|l|l|}
\hline
$N$                        & \makecell{$\mathcal{H}$-matrix} & Mem (GB) & $t_{init}$ (s) & $t_{FMVP}$ (s) & MFLOP & $t_{DMVP}$ (s)  & $RE_{MVP}$ \\ \hline
\multirow{5}{*}{64000}    & $\mathcal{H}_{hw}^{\prime}$ & 1.32               & 7.22   & 0.09           & 312.28    & 8.24 & 8.82E-06  \\ \cline{2-8} 
                         & $\mathcal{H}_{hw}^{\prime \prime}$  & \textbf{1.25}                & 15.15   & \textbf{0.08}           & \textbf{302.22}    & 8.24 & 8.78E-06  \\ \cline{2-8} 
                         & $\bkt{\mathcal{H}^2 + \mathcal{H}}_{*}$  & 1.61                & 8.52   & 0.13           & 383.94    & 8.24 & 1.07E-05  \\ \cline{2-8} 
                         & $\mathcal{H}^2_{*}$  & 1.55                & 19.72   & 0.16           & 381.14    & 8.24 & 1.05E-05  \\ \cline{2-8} 
                         & $\mathcal{H}^2_{\eta}$  & 1.77                & 19.75   & 0.18           & 435.81    & 8.24 & 8.13E-06  \\ \hline \hline
\multirow{5}{*}{512000}    & $\mathcal{H}_{hw}^{\prime}$ & 13.50               & 116.07   & 1.17           & 3195.32    & 527.61 & 9.03E-06  \\ \cline{2-8} 
                         & $\mathcal{H}_{hw}^{\prime \prime}$  & \textbf{12.03}                & 321.77   & \textbf{0.95}           & \textbf{2942.5}    & 527.61 & 8.52E-06  \\ \cline{2-8} 
                         & $\bkt{\mathcal{H}^2 + \mathcal{H}}_{*}$  & 17.91                & 136.68   & 1.48           & 4306.96    & 527.61 & 9.99E-06  \\ \cline{2-8} 
                         & $\mathcal{H}^2_{*}$  & 16.82                & 371.25   & 1.72           & 4140.91    & 527.61 & 1.02E-05  \\ \cline{2-8} 
                         & $\mathcal{H}^2_{\eta}$  & 20.64                & 346.29   & 2.42           & 5095.75    & 527.61 & 1.03E-05  \\ \hline \hline
\multirow{5}{*}{1000000}    & $\mathcal{H}_{hw}^{\prime}$ &    28.17      & 280.27    & 2.17          & 6799.32    & 2038.87 & 1.33E-05 \\ \cline{2-8} 
                         & $\mathcal{H}_{hw}^{\prime \prime}$  & \textbf{25.04}                & 756.33   & \textbf{2.02}           & \textbf{6228.76}    & 2038.87 & 1.44E-05 \\ \cline{2-8} 
                         & $\bkt{\mathcal{H}^2 + \mathcal{H}}_{*}$  & 35.98                & 334.58   & 3.87           & 8745.57    & 2038.87 & 9.76E-06  \\ \cline{2-8} 
                         & $\mathcal{H}^2_{*}$  & 32.64                & 851.42   & 3.03           & 8129.67    & 2038.87 & 1.21E-05  \\ \cline{2-8} 
                         & $\mathcal{H}^2_{\eta}$  & 38.02                & 752.54   & 5.26           & 9473.99    & 2038.87 & 8.93E-06  \\ \hline \hline
\multirow{5}{*}{2197000}    & $\mathcal{H}_{hw}^{\prime}$ & 61.02               & 670.82   & 5.21           & 14598.01    & 10062.11 & 9.22E-06 \\ \cline{2-8} 
                         & $\mathcal{H}_{hw}^{\prime \prime}$  & \textbf{51.30}               & 1960.85   & \textbf{5.02}           & \textbf{12675.9}    & 10062.11 & 9.45E-06 \\ \cline{2-8} 
                         & $\bkt{\mathcal{H}^2 + \mathcal{H}}_{*}$  & 92.22                & 810.27   & 8.28           & 22385.3    & 10062.11 & 7.34E-06  \\ \cline{2-8} 
                         & $\mathcal{H}^2_{*}$  & 83.15                & 2217.74   & 7.25           & 20637.3    & 10062.11 & 8.09E-06  \\ \cline{2-8} 
                         & $\mathcal{H}^2_{\eta}$  & 99.22                & 1956.33   & 12.52           & 24656.4    & 10062.11 & 1.21E-05  \\ \hline 
\end{tabular}
}
    \caption{Comparison between different $\mathcal{H}$-matrix formats when $\epsilon = 10^{-5}$.}
    \label{table:smvp3d_lap_2}
\end{table}

\begin{table}
      \centering
      \resizebox{\textwidth}{!}{%
      \setlength\extrarowheight{0.9pt}
\begin{tabular}{|l|l|l|l|l|l|l|l|}
\hline
$N$                        & \makecell{$\mathcal{H}$-matrix} & Mem (GB) & $t_{init}$ (s) & $t_{FMVP}$ (s) & MFLOP & $t_{DMVP}$ (s)  & $RE_{MVP}$ \\ \hline
\multirow{5}{*}{64000}    & $\mathcal{H}_{hw}^{\prime}$ & 1.99               & 31.02   & 0.17           & 469.93    & 8.24 & 6.23E-08  \\ \cline{2-8} 
                         & $\mathcal{H}_{hw}^{\prime \prime}$  & \textbf{1.87}                & 37.68   & \textbf{0.16}           & \textbf{458.75}    & 8.24 & 6.25E-08  \\ \cline{2-8} 
                         & $\bkt{\mathcal{H}^2 + \mathcal{H}}_{*}$  & 2.74                & 38.82   & 0.26           & 658.194    & 8.24 & 8.41E-08  \\ \cline{2-8} 
                         & $\mathcal{H}^2_{*}$  & 2.66                & 47.69   & 0.25           & 656.56    & 8.24 & 8.40E-08  \\ \cline{2-8} 
                         & $\mathcal{H}^2_{\eta}$  & 2.97                & 28.64   & 0.26           & 734.94    & 8.24 & 8.88E-08  \\ \hline \hline
\multirow{5}{*}{512000}    & $\mathcal{H}_{hw}^{\prime}$ & 23.34               & 673.86   & 2.20           & 5580.93    & 527.61 & 1.06E-07  \\ \cline{2-8} 
                         & $\mathcal{H}_{hw}^{\prime \prime}$  & \textbf{21.44}                & 970.29   & \textbf{1.92}           & \textbf{5294.2}    & 527.61 & 1.09E-07  \\ \cline{2-8} 
                         & $\bkt{\mathcal{H}^2 + \mathcal{H}}_{*}$  & 32.21                & 738.83   & 3.06           & 7796.51    & 527.61 & 1.48E-07  \\ \cline{2-8} 
                         & $\mathcal{H}^2_{*}$  & 30.66                & 1061.46   & 3.05           & 7599.96    & 527.61 & 1.45E-07  \\ \cline{2-8} 
                         & $\mathcal{H}^2_{\eta}$  & 35.33                & 542.03   & 3.12           & 8770.46    & 527.61 & 9.98E-08  \\ \hline \hline
\multirow{5}{*}{1000000}    & $\mathcal{H}_{hw}^{\prime}$ &    66.94      & 1711.56    & 6.64          & 16322.1    & 2038.87 & 1.63E-07 \\ \cline{2-8} 
                         & $\mathcal{H}_{hw}^{\prime \prime}$  & \textbf{62.24}                & 2523.31   & \textbf{6.18}           & \textbf{15525.8}    & 2038.87 & 1.65E-07 \\ \cline{2-8} 
                         & $\bkt{\mathcal{H}^2 + \mathcal{H}}_{*}$  & 70.77                & 1834.98   & 7.29           & 17278.6    & 2038.87 & 1.02E-07  \\ \cline{2-8} 
                         & $\mathcal{H}^2_{*}$  & 65.61                & 2660.71   & 6.98           & 16365.9    & 2038.87 & 6.40E-07  \\ \cline{2-8} 
                         & $\mathcal{H}^2_{\eta}$  & 77.98                & 1534.66   & 7.42           & 19465.5    & 2038.87 & 2.34E-07  \\ \hline \hline
\multirow{5}{*}{2197000}    & $\mathcal{H}_{hw}^{\prime}$ & 131.78               & 4191.8   & 14.51           & 31838.8    & 10062.11 & 2.35E-07 \\ \cline{2-8} 
                         & $\mathcal{H}_{hw}^{\prime \prime}$  & \textbf{116.71}               & 6682.92   & \textbf{12.73}           & \textbf{29025.5}    & 10062.11 & 2.23E-07 \\ \cline{2-8} 
                         & $\bkt{\mathcal{H}^2 + \mathcal{H}}_{*}$  & 156.17                & 4324.74   & 15.62           & 37877.4    & 10062.11 & 1.13E-07  \\ \cline{2-8} 
                         & $\mathcal{H}^2_{*}$  & 140.92                & 7012.1   & 15.02           & 35076.3    & 10062.11 & 1.15E-07  \\ \cline{2-8} 
                         & $\mathcal{H}^2_{\eta}$  & 164.58                & 3244.51   & 16.32           & 40996.1    & 10062.11 & 2.05E-07  \\ \hline 
\end{tabular}
}
    \caption{Comparison between different $\mathcal{H}$-matrix formats when $\epsilon = 10^{-7}$.}
    \label{table:smvp3d_lap_3}
\end{table}

From the above tables, we can conclude the following.
\begin{enumerate}
    \item $\mathcal{H}_{hw}^{\prime}$ and $\mathcal{H}_{hw}^{\prime \prime}$ outperform their (uniform) weak-admissibility-based counterparts, $\bkt{\mathcal{H}^2 + \mathcal{H}}_{*}$ and $\mathcal{H}^2_{*}$, respectively in both memory and MVP time, as $\mathcal{S}(\ell)>0$ in all cases (see \cref{eq:comparison_eqn}). \Cref{fig:MFLOP_comp} illustrates the dense and low-rank MFLOP of the $\mathcal{H}_{hw}^{\prime \prime}$ and $\mathcal{H}^2_{*}$ matrices for different tolerances $(\epsilon)$. It can be seen that, despite employing nested bases at all levels, $\mathcal{H}^2_{*}$ can be improved in terms of memory requirements and MVP time by using non-nested bases and avoiding the computational overhead of many translation operators at the finer levels.

    \begin{figure}
        \centering
        \subfloat[\scriptsize  $N = 512000$.]{
        \includegraphics[scale=.365]{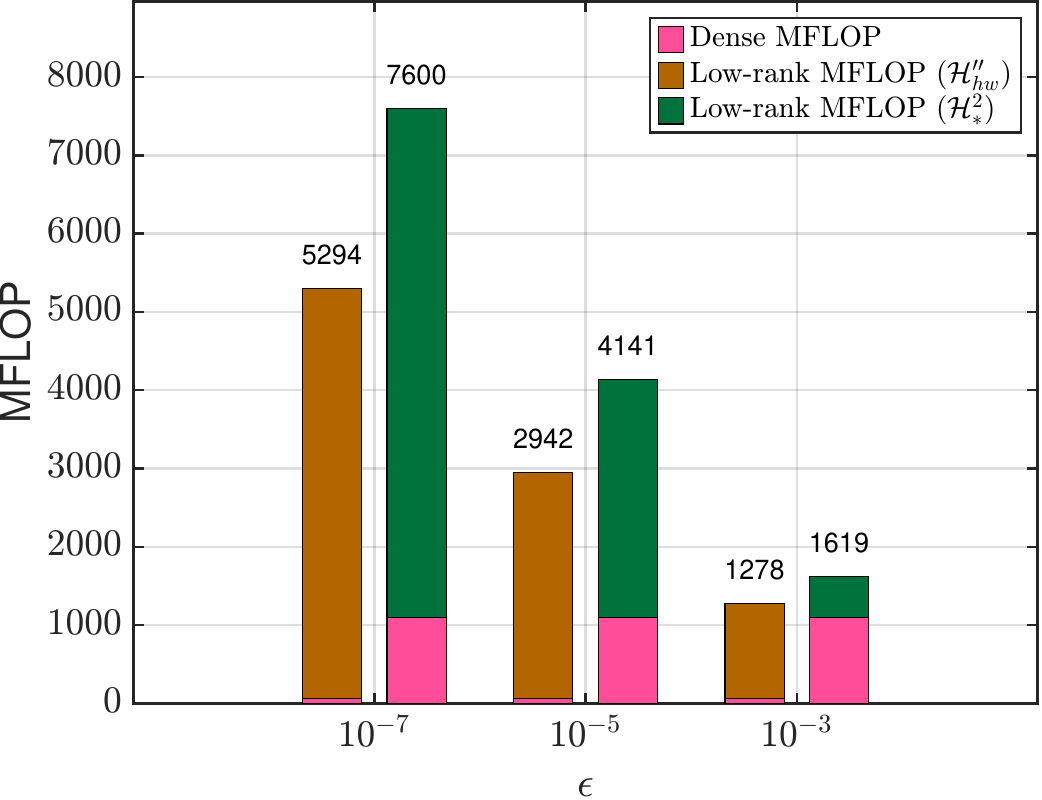}\label{fig:MFLOP1}}\space
        \subfloat[\scriptsize  $N = 2197000$.]{
        \includegraphics[scale=.365]{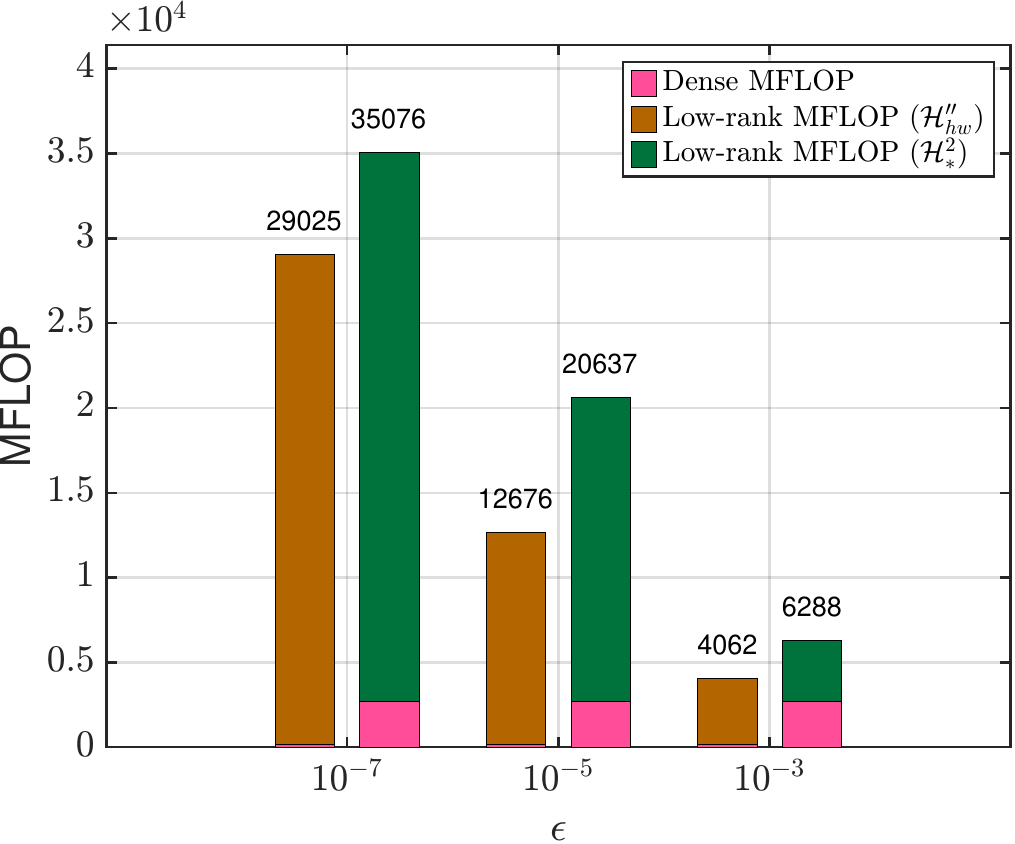}\label{fig:MFLOP2}}
        \caption{Comparison of the MFLOP of $\mathcal{H}_{hw}^{\prime\prime}$ and $\mathcal{H}^{2}_*$. For $\mathcal{H}_{hw}^{\prime\prime}$, only the diagonal blocks are dense, whereas both neighboring and diagonal blocks are dense for $\mathcal{H}^{2}_*$.}
        \label{fig:MFLOP_comp}
    \end{figure}

    \item  At finer levels, the construction of $\mathcal{H}_{hw}^{\prime}$ and $\mathcal{H}_{hw}^{\prime \prime}$ is similar. However, $\mathcal{H}_{hw}^{\prime \prime}$ achieves better memory efficiency and MVP performance due to the use of nested bases for vertex-sharing interactions at coarser levels. This can also be observed in other experiments.

    \item As shown in \cite{nHODLRdD}, $\bkt{\mathcal{H}^2 + \mathcal{H}}_{*}$ and $\mathcal{H}^2_{*}$ outperform the standard admissibility based $\mathcal{H}^2_{\eta}$ in terms of memory and MVP time. Arranging the five $\mathcal{H}$-matrix formats (\Cref{tab:H_mat_summ}) in ascending order of memory usage and MVP time yields

    $$ \mathcal{H}_{hw}^{\prime \prime} < \mathcal{H}_{hw}^{\prime} < \mathcal{H}^2_{*} < \bkt{\mathcal{H}^2 + \mathcal{H}}_{*} < \mathcal{H}^2_{\eta}.$$

   This demonstrates the superior efficiency of the proposed hybrid weak admissibility condition (\cref{eq:hybrid_admiss}) compared with both (uniform) weak admissibility (\cref{eq:weak_admis}) and standard admissibility conditions (\cref{eq:H_admiss}).

    \end{enumerate}

\subsection{Fast principal component analysis}
Gaussian process (GP) is a probabilistic model used in machine learning and statistics to predict outcomes. However, when working with high-dimensional data, GP can become computationally expensive because matrix operations scale poorly with the size and dimensionality of the data. To address this issue, dimension reduction techniques such as principal component analysis are essential. Principal component analysis (PCA) helps by projecting the input data onto a lower-dimensional subspace, retaining the most important information while reducing the number of dimensions. In GP, the covariance matrix is a key component representing the relationships between data points based on the kernel. 

We use randomized Nystr\"om approximation \cite{martinsson2020randomized} to obtain truncated eigenvalue decomposition of the covariance matrix. In the naive randomized Nystr\"om algorithm, computing the top-$k$ eigenvalues and their corresponding eigenvectors requires a computational cost of $\mathcal{O}\bkt{kN^2}$, where the covariance matrix $K\in \mathbb{R}^{N \times N}$. The most computationally expensive part of the algorithm is performing matrix-matrix products, such as $K \Omega$ or $K Q$. In the fast PCA (see \Cref{sec:fpca}), we replace $K$ with its hierarchical representations ($\Tilde{K}$), enabling fast hierarchical matrix-vector products (refer to lines 4 and 5 in \Cref{algo:range_finder} and \Cref{algo:nystrom_approx}, respectively). Hence, it reduces the computational cost to $\mathcal{O}\bkt{kN \log^{\alpha}(N)} $, where $\alpha \geq 0$. 

In this experiment, we consider the Matérn-$3/2$ covariance kernel and the Gaussian kernel $(e^{-r^2})$ with $N$ randomly generated points in the cube $[-1,1]^3$. We set $k=100$ and $c=10$ in \Cref{algo:nystrom_approx} and compute top-$100$ eigenpairs of the covariance matrices. In total, $220$ MVPs are performed. For the fast PCA method, we set the tolerance $\epsilon = 10^{-4}$ to construct $\Tilde{K}$. In the direct PCA algorithm, MVPs are carried out using the \href{https://eigen.tuxfamily.org}{Eigen} library. 

We present the performance of the fast PCA over direct PCA in \Cref{tab:pca_ex1}. In addition, we report the memory and relative error ($RE_{PCA}$) $\magn{\Lambda_{\Tilde{K}} - \Lambda_{K}}_2/\magn{\Lambda_{K}}_2$, where $\Lambda_{\Tilde{K}}$ and $\Lambda_{K}$ represent the eigenvalues computed using fast PCA and direct PCA, respectively. As observed from \Cref{tab:pca_ex1}, the memory and time of fast PCA scales quasi-linearly with $N$, whereas the direct PCA demonstrates quadratic scaling.

\begin{table}[H]
\centering
\resizebox{\textwidth}{!}{%
\setlength\extrarowheight{0.9pt}
\begin{tabular}{|l|l|c|llll|l|}
\hline
\multirow{2}{*}{$N$} &
\multirow{2}{*}{\makecell{$\mathcal{H}$-matrix}} &
\multirow{2}{*}{\makecell{Mem(GB)}} &
\multicolumn{4}{c|}{Time (s)} &
\multirow{2}{*}{$RE_{PCA}$} \\ \cline{4-7}

&
&
&
$t_{init}$ &
$t_{FPCA}$ &
$t_{total}$ &
$t_{DPCA}$ &
\\ \hline

\multirow{2}{*}{1000}
& $\mathcal{H}_{hw}^{\prime}$
& 0.002
& 0.005 & 0.052 & 0.057 & 1.58 & 2.89E-06 \\ \cline{2-8}

& $\mathcal{H}_{hw}^{\prime \prime}$
& \textbf{0.002}
& 0.005 & 0.051 & \textbf{0.056} & 1.58 & 3.04E-06 \\ \hline\hline

\multirow{2}{*}{8000}
& $\mathcal{H}_{hw}^{\prime}$
& 0.026
& 0.062 & 0.607 & 0.669 & 94.11 & 2.67E-06 \\ \cline{2-8}

& $\mathcal{H}_{hw}^{\prime \prime}$
& \textbf{0.025}
& 0.066 & 0.601  & \textbf{0.667} & 94.11 & 3.15E-06 \\ \hline \hline

\multirow{2}{*}{64000}
& $\mathcal{H}_{hw}^{\prime}$
& 0.216
& 1.781 & 6.101 & 7.882 & 5592.79 & 8.07E-06 \\ \cline{2-8}

& $\mathcal{H}_{hw}^{\prime \prime}$
& \textbf{0.193}
& 1.924 & 5.844 & \textbf{7.768} & 5592.79 & 8.39E-06 \\ \hline\hline

\multirow{2}{*}{512000}
& $\mathcal{H}_{hw}^{\prime}$
& 1.668
& 20.928 & 53.399 & 74.327 & -- & -- \\ \cline{2-8}

& $\mathcal{H}_{hw}^{\prime \prime}$
& \textbf{1.281}
& 23.044 & 47.058 & \textbf{70.102} & -- & -- \\ \hline

\end{tabular}
}
\caption{Comparison between $\mathcal{H}$-matrix accelerated PCA $(FPCA)$ and direct PCA $(DPCA)$ for Matérn-$3/2$ covariance kernel. $t_{DPCA}$ scales as $\mathcal{O}(N^2)$.}
\label{tab:pca_ex1}
\end{table}

\begin{table}[H]
\centering
\resizebox{\textwidth}{!}{%
\setlength\extrarowheight{0.9pt}
\begin{tabular}{|l|l|c|llll|l|}
\hline
\multirow{2}{*}{$N$} &
\multirow{2}{*}{\makecell{$\mathcal{H}$-matrix}} &
\multirow{2}{*}{\makecell{Mem(GB)}} &
\multicolumn{4}{c|}{Time (s)} &
\multirow{2}{*}{$RE_{PCA}$} \\ \cline{4-7}

&
&
&
$t_{init}$ &
$t_{FPCA}$ &
$t_{total}$ &
$t_{DPCA}$ &
\\ \hline

\multirow{2}{*}{1000}
& $\mathcal{H}_{hw}^{\prime}$
& 0.003
& 0.007 & 0.058 & 0.065 & 1.52 & 7.52E-06 \\ \cline{2-8}

& $\mathcal{H}_{hw}^{\prime \prime}$
& \textbf{0.002}
& 0.006 & 0.057 & \textbf{0.063} & 1.52 & 7.59E-06 \\ \hline\hline

\multirow{2}{*}{8000}
& $\mathcal{H}_{hw}^{\prime}$
& 0.036
& 0.083 & 0.813 & 0.896 & 92.17 & 4.90E-06 \\ \cline{2-8}

& $\mathcal{H}_{hw}^{\prime \prime}$
& \textbf{0.034}
& 0.083 & 0.801 & \textbf{0.884} & 92.17 & 5.77E-06 \\ \hline\hline

\multirow{2}{*}{64000}
& $\mathcal{H}_{hw}^{\prime}$
& 0.310
& 3.054 & 7.907 & 10.961 & 5158.39 & 2.76E-05 \\ \cline{2-8}

& $\mathcal{H}_{hw}^{\prime \prime}$
& \textbf{0.262}
& 3.378 & 7.048 & \textbf{10.426} & 5158.39 & 1.78E-05 \\ \hline\hline

\multirow{2}{*}{512000}
& $\mathcal{H}_{hw}^{\prime}$
& 2.353
& 28.462 & 65.572 & 94.034 & -- & -- \\ \cline{2-8}

& $\mathcal{H}_{hw}^{\prime \prime}$
& \textbf{1.660}
& 33.917 & 54.452 & \textbf{88.369} & -- & -- \\ \hline

\end{tabular}
}
\caption{Comparison between $\mathcal{H}$-matrix accelerated PCA $(FPCA)$ and direct PCA $(DPCA)$ for Gaussian kernel. $t_{DPCA}$ scales as $\mathcal{O}(N^2)$.}
\label{tab:pca_ex2}
\end{table}

\subsection{Fast kernel ridge regression}
We discuss a fast iterative solution for regularized linear systems \cref{eq:num_result_krr}, which frequently arise in many applications, such as kernel ridge regression (a supervised learning technique) \cite{alaoui2015fast}.

\begin{equation} \label{eq:num_result_krr}
    \bkt{K + \mu I} \xb = \pmb{b},
\end{equation}
where $K \in \mathbb{R}^{N \times N}$ is a symmetric positive definite (SPD) kernel matrix and regularization parameter $\mu \geq 0$. 

\cref{eq:krr} is typically ill-conditioned, leading to slow convergence for the iterative solver. Hence, a preconditioner for the iterative solver is essential. We use randomized Nystr\"om preconditioner \cite{frangella2023randomized} here. Thus, \cref{eq:num_result_krr} is solved in the following steps.

\begin{enumerate}
    \item We construct the \textbf{fast} randomized Nystr\"om preconditioner $P$ as described in \Cref{sec:nystrom_precon}.
    \item Once we have the preconditioner $P$, we use it in the \textbf{fast} preconditioned conjugate gradient (PCG) method. This is outlined in \Cref{sec:pcg}.
\end{enumerate}

In both steps, we employ the hierarchical representation of $K$, i.e., $\Tilde{K}$, to accelerate the associated MVPs.

For a numerical experiment, we consider the Matérn-$1/2$ kernel $(\exp{(-r)})$ with $\mu = 0.0001$. Synthetic datasets are generated by randomly distributing points, $\pmb{x}_i$, in the cube $[-1,1]^3$. We use $N$ points for training and $m=N/5$ additional points for testing. The corresponding output values are denoted by $\pmb{b}_i \in \Rb$. The hierarchical representation $\Tilde{K}$ is formed by setting tolerance $\epsilon = 10^{-5}$. $\Tilde{K}$ is used to construct the fast randomized Nystr\"om preconditioner $P$ and to accelerate the MVP in CG. While constructing the fast randomized Nystr\"om preconditioner, we set rank $k=250$. For CG, the maximum number of iterations is $500$ with stopping tolerance $\epsilon_{CG} = 10^{-5}$.

The performance of fast Nystr\"om PCG is presented in \Cref{tab:npcg_ex1}. We report the root mean squared error (RMSE) between the predicted and true values. $t_{FNPCG}$ represents the time taken to construct the fast Nyström preconditioner and perform the fast PCG. We also benchmark the solution of \cref{eq:num_result_krr} using the diagonal preconditioned conjugate gradient (DPCG) and unpreconditioned conjugate gradient (CG) methods and report the number of iterations for each. It can be noted from \Cref{tab:npcg_ex1} that for $N>1000$, both DPCG and CG fail to converge within $500$ iterations, whereas the fast Nystr\"om PCG (FNPCG) method converges.

\begin{table}[H]
\centering
\resizebox{\textwidth}{!}{%
\setlength\extrarowheight{0.9pt}
\begin{tabular}{|c|c|c|ccc|ccc|}
\hline
\multirow{2}{*}{$N$} &
\multirow{2}{*}{\makecell{$\mathcal{H}$-matrix}} &
\multirow{2}{*}{RMSE} &
\multicolumn{3}{c|}{Time (s)} &
\multicolumn{3}{c|}{\# Iter} \\ \cline{4-9}

& & & $t_{init}$ & $t_{FNPCG}$ & $t_{total}$
& FNPCG & DPCG & CG \\ \hline \hline

\multirow{2}{*}{1000}
& $\mathcal{H}_{hw}^{\prime}$
& 0.2 \% & 0.008 & 0.15 & 0.158 & 9 & 112 & 156 \\ \cline{2-9}

& $\mathcal{H}_{hw}^{\prime \prime}$
& 0.2 \% & 0.008 & 0.13 & \textbf{0.138} & 10 & 112 & 156 \\ \hline \hline

\multirow{2}{*}{8000}
& $\mathcal{H}_{hw}^{\prime}$
& 0.9 \% & 0.28 & 2.78 & 3.06 & 65 & - & - \\ \cline{2-9}

& $\mathcal{H}_{hw}^{\prime \prime}$
& 0.9 \% & 0.28 & 2.36 & \textbf{2.64} & 67 & - & - \\ \hline \hline

\multirow{2}{*}{64000}
& $\mathcal{H}_{hw}^{\prime}$
& 0.8 \% & 7.48 & 56.31 & 63.79 & 312 & - & - \\ \cline{2-9}

& $\mathcal{H}_{hw}^{\prime \prime}$
& 0.8 \% & 8.08 & 52.71 & \textbf{60.79} & 316 & - & - \\ \hline 

% \multirow{2}{*}{512000}
% & $\mathcal{H}_{hw}^{\prime}$
% & 3.2 & 112.57 & 551.12 & 663.69 & 48 & - & - \\ \cline{2-9}

% & $\mathcal{H}_{hw}^{\prime \prime}$
% & 3.2 & 140.31 & 592.58 & 732.89 & 60 & - & - \\ \hline

\end{tabular}
}
\caption{Numerical results for the synthetic dataset with Matérn-$1/2$ kernel. 
We report the number of iterations for fast Nyström PCG (FNPCG), diagonal PCG (DPCG), 
and unpreconditioned CG (CG). The symbol ``$-$'' indicates that the iterative solver 
fails to converge within $500$ iterations.}
\label{tab:npcg_ex1}
\end{table}

\section{Conclusion} \label{sec:conclusion}
We have presented a new hybrid weak admissibility condition, under which far-field and vertex-sharing clusters are admissible at coarser levels, and adjacent clusters are considered admissible at finer levels. Based on this admissibility condition, we have proposed two new hierarchical matrix formats, $\mathcal{H}_{hw}^{\prime}$ and $\mathcal{H}_{hw}^{\prime \prime}$. For the $\mathcal{H}_{hw}^{\prime}$ format, we have employed nested bases for far-field interactions and non-nested bases for vertex-sharing interactions. In contrast, for the $\mathcal{H}_{hw}^{\prime \prime}$ format, we have employed nested bases for both far-field and vertex-sharing interactions. In both formats, we have used non-nested bases at the finer levels. Numerical experiments have demonstrated that the proposed hybrid hierarchical matrix formats outperform their counterparts based on uniform weak admissibility, as well as the standard $\mathcal{H}^2$-matrix format, thereby illustrating the effectiveness of the proposed hybrid weak admissibility condition. We have further applied the proposed hierarchical matrix formats to the acceleration of principal component analysis in Gaussian processes and kernel ridge regression problems. The \texttt{C++} implementation of the proposed hierarchical matrix algorithms is publicly available at \url{https://github.com/riteshkhan/Hhweak}.

\section*{Acknowledgment}
R.K. acknowledges support from the European Union (ERC, inEXASCALE, Grant No. 101075632) and the Charles University Research Centre program (No. UNCE/24/SCI/005).

\bibliographystyle{siamplain}
\bibliography{refs}

\begin{thebibliography}{10}

\bibitem{alaoui2015fast}
{\sc A.~Alaoui and M.~W. Mahoney}, {\em Fast randomized kernel ridge regression with statistical guarantees}, Advances in neural information processing systems, 28 (2015).

\bibitem{ambartsumyan2020hierarchical}
{\sc I.~Ambartsumyan, W.~Boukaram, T.~Bui-Thanh, O.~Ghattas, D.~Keyes, G.~Stadler, G.~Turkiyyah, and S.~Zampini}, {\em Hierarchical matrix approximations of hessians arising in inverse problems governed by pdes}, SIAM Journal on Scientific Computing, 42 (2020), pp.~A3397--A3426.

\bibitem{ambikasaran2013mathcal}
{\sc S.~Ambikasaran and E.~Darve}, {\em An $\mathcal{O}\left({N} \log {N}\right)$ fast direct solver for partial hierarchically semi-separable matrices}, Journal of Scientific Computing, 57 (2013), pp.~477--501.

\bibitem{ambikasaran2025hybrid}
{\sc S.~Ambikasaran, R.~Khan, J.~Tausch, and S.~Wang}, {\em A hybrid interpolation aca accelerated method for parabolic boundary integral operators}, SIAM Journal on Scientific Computing, 47 (2025), pp.~A1507--A1526.

\bibitem{ambikasaran2019hodlrlib}
{\sc S.~Ambikasaran, K.~R. Singh, and S.~S. Sankaran}, {\em {HODLRlib}: {A} library for hierarchical matrices}, Journal of Open Source Software, 4 (2019), p.~1167.

\bibitem{amestoy2017complexity}
{\sc P.~Amestoy, A.~Buttari, J.-Y. L'Excellent, and T.~Mary}, {\em On the complexity of the block low-rank multifrontal factorization}, SIAM Journal on Scientific Computing, 39 (2017), pp.~A1710--A1740.

\bibitem{aca}
{\sc M.~Bebendorf and S.~Rjasanow}, {\em Adaptive low-rank approximation of collocation matrices}, Computing, 70 (2003), pp.~1--24.

\bibitem{bebendorf2012constructing}
{\sc M.~Bebendorf and R.~Venn}, {\em Constructing nested bases approximations from the entries of non-local operators}, Numerische Mathematik, 121 (2012), pp.~609--635.

\bibitem{borm2009construction}
{\sc S.~B{\"o}rm}, {\em Construction of data-sparse $\mathcal{H}^2$-matrices by hierarchical compression}, SIAM Journal on Scientific Computing, 31 (2009), pp.~1820--1839.

\bibitem{borm2003hierarchical}
{\sc S.~B{\"o}rm, L.~Grasedyck, and W.~Hackbusch}, {\em Hierarchical matrices}, Lecture notes, 21 (2003), p.~2003.

\bibitem{chandrasekaran2007fast}
{\sc S.~Chandrasekaran, P.~Dewilde, M.~Gu, W.~Lyons, and T.~Pals}, {\em A fast solver for {HSS} representations via sparse matrices}, SIAM Journal on Matrix Analysis and Applications, 29 (2007), pp.~67--81.

\bibitem{frangella2023randomized}
{\sc Z.~Frangella, J.~A. Tropp, and M.~Udell}, {\em Randomized nystr{\"o}m preconditioning}, SIAM Journal on Matrix Analysis and Applications, 44 (2023), pp.~718--752.

\bibitem{gillman2012direct}
{\sc A.~Gillman, P.~M. Young, and P.-G. Martinsson}, {\em A direct solver with $\mathcal{O} ({N})$ complexity for integral equations on one-dimensional domains}, Frontiers of Mathematics in China, 7 (2012), pp.~217--247.

\bibitem{grasedyck2003construction}
{\sc L.~Grasedyck and W.~Hackbusch}, {\em Construction and arithmetics of h-matrices}, Computing, 70 (2003), pp.~295--334.

\bibitem{gray2000}
{\sc A.~Gray and A.~Moore}, {\em N-body'problems in statistical learning}, Advances in neural information processing systems, 13 (2000).

\bibitem{greengard1987fast}
{\sc L.~Greengard and V.~Rokhlin}, {\em A fast algorithm for particle simulations}, Journal of computational physics, 73 (1987), pp.~325--348.

\bibitem{gujjula2022new}
{\sc V.~Gujjula and S.~Ambikasaran}, {\em A new nested cross approximation}, arXiv preprint arXiv:2203.14832,  (2022).

\bibitem{hackbusch1999sparse}
{\sc W.~Hackbusch}, {\em A sparse matrix arithmetic based on $\mathcal{H}$-matrices. part i: Introduction to $\mathcal{H}$-matrices}, Computing, 62 (1999), pp.~89--108.

\bibitem{hmatrix_book}
{\sc W.~Hackbusch}, {\em Hierarchical Matrices: Algorithms and Analysis}, vol.~49, Springer, 12 2015, \url{https://doi.org/10.1007/978-3-662-47324-5}.

\bibitem{hackbusch2000sparse}
{\sc W.~Hackbusch and B.~N. Khoromskij}, {\em A sparse $\mathcal{H}$-matrix arithmetic. part ii: Application to multi-dimensional problems}, Computing, 64 (2000), pp.~21--47.

\bibitem{weak_hackbusch2004}
{\sc W.~Hackbusch, B.~N. Khoromskij, and R.~Kriemann}, {\em Hierarchical matrices based on a weak admissibility criterion}, Computing, 73 (2004), pp.~207--243.

\bibitem{halko2011finding}
{\sc N.~Halko, P.~G. Martinsson, and J.~A. Tropp}, {\em Finding structure with randomness: Probabilistic algorithms for constructing approximate matrix decompositions}, SIAM review, 53 (2011), pp.~217--288.

\bibitem{harlim2020kernel}
{\sc J.~Harlim, D.~Sanz-Alonso, and R.~Yang}, {\em Kernel methods for bayesian elliptic inverse problems on manifolds}, SIAM/ASA Journal on Uncertainty Quantification, 8 (2020), pp.~1414--1445.

\bibitem{ho2013hierarchical}
{\sc K.~L. Ho and L.~Ying}, {\em Hierarchical interpolative factorization for elliptic operators: integral equations}, arXiv preprint arXiv:1307.2666,  (2013).

\bibitem{kandappan2022hodlr2d}
{\sc V.~A. Kandappan, V.~Gujjula, and S.~Ambikasaran}, {\em {HODLR2D}: A new class of hierarchical matrices}, arXiv preprint arXiv:2204.05536,  (2022).

\bibitem{nHODLRdD}
{\sc R.~Khan and S.~Ambikasaran}, {\em New algebraic fast algorithms for ${N}$-body problems in two and three dimensions}, Communications in Computational Physics, 37 (2025), pp.~1157--1226.

\bibitem{khan2022numerical}
{\sc R.~Khan, V.~A. Kandappan, and S.~Ambikasaran}, {\em {HODLRdD}: A new black-box fast algorithm for ${N}$-body problems in d-dimensions with guaranteed error bounds: Applications to integral equations and support vector machines}, Journal of Computational Physics, 501 (2024), p.~112786, \url{https://doi.org/10.1016/j.jcp.2024.112786}.

\bibitem{martinsson2020randomized}
{\sc P.-G. Martinsson and J.~A. Tropp}, {\em Randomized numerical linear algebra: Foundations and algorithms}, Acta Numerica, 29 (2020), pp.~403--572.

\bibitem{massei2022hierarchical}
{\sc S.~Massei, L.~Robol, and D.~Kressner}, {\em Hierarchical adaptive low-rank format with applications to discretized partial differential equations}, Numerical Linear Algebra with Applications, 29 (2022), p.~e2448.

\bibitem{gp_book}
{\sc C.~E. Rasmussen}, {\em Gaussian Processes in Machine Learning}, Springer Berlin Heidelberg, Berlin, Heidelberg, 2004, pp.~63--71, \url{https://doi.org/10.1007/978-3-540-28650-9_4}, \url{https://doi.org/10.1007/978-3-540-28650-9_4}.

\bibitem{xia2010fast}
{\sc J.~Xia, S.~Chandrasekaran, M.~Gu, and X.~S. Li}, {\em Fast algorithms for hierarchically semiseparable matrices}, Numerical Linear Algebra with Applications, 17 (2010), pp.~953--976.

\bibitem{ker_den}
{\sc C.~Yang, R.~Duraiswami, N.~A. Gumerov, and L.~Davis}, {\em Improved fast gauss transform and efficient kernel density estimation}, in Computer Vision, IEEE International Conference on, vol.~2, IEEE Computer Society, 2003, pp.~464--464.

\bibitem{zhao2019fast}
{\sc Y.~Zhao, D.~Jiao, and J.~Mao}, {\em Fast nested cross approximation algorithm for solving large-scale electromagnetic problems}, IEEE Transactions on Microwave Theory and Techniques, 67 (2019), pp.~3271--3283.

\end{thebibliography}

\newpage
\appendix
\section{Fast principal component analysis} \label{sec:fpca}
The randomized Nystr\"om approximation algorithm can be found in \cite{halko2011finding, martinsson2020randomized}. In the original algorithm, the most computationally expensive part is matrix-matrix products, such as $K \Omega$ or $K Q$. We accelerate these operations by replacing the covariance matrix $K$ with its hierarchical representations $\Tilde{K}$, leading to fast principal component analysis (see \Cref{algo:range_finder} and \Cref{algo:nystrom_approx}). 

\begin{algorithm}[H]
\small
\caption{Fast randomized range finder}
\label{algo:range_finder}
\begin{algorithmic}[1]
    \State \textbf{Input:} Hierarchical representation of $K\in \mathbb{R}^{N \times N}$, denoted as $\Tilde{K}$ and $l$
    \State \textbf{Output:}  $Q \in \mathbb{R}^{N \times l}$
\Function{Fast Randomized Range Finder}{$\Tilde{K}, l$}
    \State  $\Omega = \texttt{randn}\bkt{N,l}$
    \State  $Y = \Tilde{K} \Omega$  \Comment{\textbf{We perform $l$ fast MVPs}}
    \State $Q = \texttt{qr} \bkt{Y,0}$ \Comment{Thin QR factorization}
\EndFunction
\end{algorithmic}
\end{algorithm}

\begin{algorithm}[H]
\small
\caption{Fast Nystr\"om approximation of $K \in \mathbb{R}^{N \times N}$, i.e., $K_{nys} = U \hat{\Lambda} U^{T}$}
\label{algo:nystrom_approx}
\begin{algorithmic}[1]
    \State \textbf{Input:} Hierarchical representation of SPD matrix $K$, $\Tilde{K}$ and rank $l = k + c$ 
    \State \textbf{Output:} Rank-$k$ approximation $K_{nys} = U \hat{\Lambda} U^{T}$, where  $U \in \mathbb{R}^{N \times k}, \hat{\Lambda} \in \mathbb{R}^{k \times k}$
\Function{Fast Nystr\"om Approximation}{$\Tilde{K}, l$}
    \State $Q = \text{Fast Randomized Range Finder}\bkt{\Tilde{K}, l}$ \Comment{From \Cref{algo:range_finder}}
    \State $B_1 = \Tilde{K} Q$  \Comment{\textbf{We perform $l$ fast MVPs}}
    \State $\nu  = \texttt{eps} \bkt{\texttt{norm}(B_1)}$  \Comment{Compute shift}
    \State $B_{\nu} = \bkt{B_1 + \nu Q}$
    \State $B_2 = Q^{*} B_{\nu}$
    \State $C = \texttt{chol}(B_2)$    \Comment{Cholesky factorization}
    \State $F = B_{\nu} / C$           \Comment{Solve using a triangular solve}
    \State $ [U, \Sigma, \sim] = \texttt{svd}\bkt{F,0}$ \Comment{Thin SVD}
    \State $\Lambda = \max\{0, \Sigma^{2} -\nu I \}$. \Comment{Remove shift}
    \State $U = U(:,1:k)$ and $\hat{\Lambda} = \Lambda(1:k, 1:k)$
\EndFunction
\end{algorithmic}
\end{algorithm}

% \begin{algorithm}[H]
% \caption{Fast eigenvalue decomposition of $K$, i.e., $K \approx U \Lambda U^{*}$}
% \label{eig_decom}
% \begin{algorithmic}[1]
%     \State \textbf{Input:} $K \in \mathbb{R}^{N \times N}, Q$ ($Q$ is obtained from \Cref{range_finder})
%     \State \textbf{Output:}  $U, \Lambda$
% \Procedure{Fast Nystr\"om Method}{}
%     \State $B_1 = K Q$  \Comment{\textbf{We perform Fast MVPs}}
%     \State $B_2 = Q^{*} B_1$
%     \State $C = \texttt{chol}(B_2)$    \Comment{Cholesky factorization}
%     \State $F = B_1 / C$ \Comment{Solve using a triangular solve}
%     \State $ [U, \Sigma, \sim] = \texttt{svd}\bkt{F,0}$ \Comment{Thin SVD}
%     \State $\Lambda = \Sigma^{2}$
% \EndProcedure
% \end{algorithmic}
% \end{algorithm}

\section{Fast kernel ridge regression}
Regularized linear systems, like \cref{eq:krr}, arise in many applications, including kernel ridge regression.

\begin{equation} \label{eq:krr}
    \bkt{K + \mu I} \xb = \pmb{b}
\end{equation}
where $K \in \mathbb{R}^{N \times N}$ is a symmetric positive definite (SPD) kernel matrix and regularization parameter $\mu \geq 0$.

 We leverage the hierarchical representation of $K$, i.e., $\Tilde{K}$, to construct the Nystr\"om preconditioner (see \Cref{sec:nystrom_precon}) and accelerate MVPs in the conjugate gradient algorithm (see \Cref{sec:pcg}).

\subsection{Fast randomized Nystr\"om preconditioner} \label{sec:nystrom_precon}
The randomized Nystr\"om preconditioner and its inverse are given as follows:

\begin{equation}
    P = \dfrac{1}{\hat{\lambda_k} + \mu} U \bkt{\hat{\Lambda} + \mu I} U^{T} + \bkt{I - U U^{T}}  
\end{equation}
  
\begin{equation}
    P^{-1} = \bkt{\hat{\lambda_k} + \mu} U \bkt{\hat{\Lambda} + \mu I}^{-1} U^{T} + \bkt{I - U U^{T}},  
\end{equation}

where $K_{nys} = U \hat{\Lambda} U^{T}$ and $\hat{\lambda_k}$ is the $k^{th}$ eigenvalue.

The naive cost to construct the Nystr\"om preconditioner $P$ scales as $\mathcal{O}\bkt{kN^2}$. However, the factors associated with the Nyström preconditioner, $U$ and $\hat{\Lambda}$, are obtained using the fast Nystr\"om Approximation (see \Cref{algo:nystrom_approx}), leading to computational cost $\mathcal{O}\bkt{kN \log^{\alpha}{(N)}}$, $\alpha \geq 0$. 

\subsection{Fast PCG} \label{sec:pcg}
The fast randomized Nystr\"om preconditioner is used as a left-preconditioner in the CG algorithm. The MVPs in PCG are accelerated using the hierarchical representation (see \Cref{algo:cg}), resulting in fast PCG.

\begin{algorithm}[H]
\small
\caption{Fast PCG} \label{algo:pcg}
\label{algo:cg}
\begin{algorithmic}[1]
    \State \textbf{Input:} Hierarchical representation of the matrix $K$, denoted as $\Tilde{K}$, regularization parameter $\mu$, right hand side $\pmb{b}$, initial guess $\pmb{x_0}$, CG tolerance $\epsilon_{CG}$ and $P^{-1}$. 
    \State \textbf{Output:}  Approximate solution $\hat{x}$
\Function{PCG}{$\Tilde{K}$, $\mu$, $\pmb{b}$, $\pmb{x_0}$, $\epsilon_{CG}$, $P^{-1}$}
    \State $\pmb{r_0} = \pmb{b} - \mu \pmb{x_0} - \Tilde{K} \pmb{x}_0$       \Comment{\textbf{We perform fast MVP}}
    \State $\pmb{z_0} = P^{-1} \pmb{r_0}$                       \Comment{Apply Nystr\"om preconditioner}
    \State $\pmb{p_0} = \pmb{z_0}$
    \While{$\magn{\pmb{r}}_2 > \epsilon_{CG}$}
    \State $\pmb{v} = \Tilde{K} \pmb{p_0} + \mu \pmb{p_0}$            \Comment{\textbf{We perform fast MVP}}
    \State $\alpha = \dfrac{\pmb{r_0}^{T} \pmb{z_0}}{\pmb{p_0}^{T}\pmb{v_0}}$
    \State $\pmb{x} = \pmb{x_0} + \alpha \pmb{p_0}$
    \State $\pmb{r} = \pmb{r_0} - \alpha \pmb{v}$
    \State $\pmb{z} = P^{-1} \pmb{r}$                           \Comment{Apply Nystr\"om preconditioner}
    \State $\beta = \dfrac{\pmb{r}^{T} \pmb{z}}{\pmb{r_0}^{T}\pmb{z_0}}$
    \State $\pmb{x_0} \leftarrow  \pmb{x}, \pmb{r_0} \leftarrow  \pmb{r}, \pmb{p_0} \leftarrow  (\pmb{z} + \beta \pmb{p_0}), \pmb{z_0} \leftarrow \pmb{z}$ \Comment{Update the values}
    \EndWhile
\EndFunction
\end{algorithmic}
\end{algorithm}

The overall computational cost for constructing the fast Nystr\"om preconditioner (\Cref{sec:nystrom_precon}) and executing the fast PCG (\Cref{sec:pcg}) scales quasi-linearly with $N$.

\end{document}